\documentclass[preprint,3p,times]{elsarticle}

\usepackage{amssymb}
\usepackage{amsmath,amsfonts}
\usepackage{graphicx}
\graphicspath{{figs/}}
\usepackage[framemethod=tikz]{mdframed}
\usepackage{booktabs}
\usepackage{subcaption}
\usepackage{xfrac}
\usepackage{float}
\usepackage{color}
\usepackage{tabularx}
\newcolumntype{Y}{>{\centering\arraybackslash}X}
\usepackage[figuresright]{rotating}
\usepackage[section]{placeins} 

\makeatletter
\def\ps@pprintTitle{%
  \let\@oddhead\@empty
  \let\@evenhead\@empty
  \let\@oddfoot\@empty
  \let\@evenfoot\@oddfoot}
\makeatother

\newcommand{\code}[1]{\texttt{#1}}
\newcommand{\grad}[0]{\nabla}
\newcommand{\diverg}[0]{\grad \cdot}

\newcommand{\pdv}[2]{{\frac{\partial{#1}}{\partial{#2}}}}
\newcommand{\bv}[1]{{\boldsymbol{#1}}}
\newcommand{\grins}[0]{\code{GRINS} }
\newcommand{\libmesh}[0]{\code{libMesh} }

\begin{document}

\begin{frontmatter}

\title{A Volumetrically Stabilized Mixed Formulation of the Finite Element
Immersed Boundary Method for Fluid--Structure Interaction with Fully
Incompressible Hyperelastic Solids}

\author[ub]{Abhishek Mishra\corref{cor1}\fnref{amd}}
\ead{amishra1995@gmail.com}
\author[ub]{Paul T. Bauman\fnref{amd}}
\cortext[cor1]{Corresponding author.}
\address[ub]{Department of Mechanical and Aerospace Engineering,
University at Buffalo, Buffalo, NY 14260, USA}
\fntext[amd]{Current affiliation: Advanced Micro Devices, Inc.}

\begin{abstract}
The finite element immersed boundary method (FE-IBM) with a distributed
Lagrange multiplier is an attractive framework for fluid--structure interaction
(FSI) because it couples an Eulerian description of an incompressible fluid to a
Lagrangian description of an immersed solid on two independent, non-conforming
meshes, avoiding the costly remeshing required by body-fitted arbitrary
Lagrangian--Eulerian methods. When the immersed solid is modelled as a fully
incompressible hyperelastic material, however, a direct finite element
discretization of the deviatoric solid stress fails to enforce the Lagrangian
incompressibility constraint, and the computed structure exhibits spurious
volumetric instabilities and locking. In this work we present a
\emph{mixed formulation} of the distributed-Lagrange-multiplier FE-IBM that
restores volumetric stability. Following the theory of nearly incompressible
hyperelasticity, we augment the solid stress with a volumetric contribution
derived from a dilatational strain energy and introduce an additional solid
pressure field, enforced weakly, that plays the role of the Lagrange multiplier
for the Lagrangian incompressibility constraint~$J=1$. The resulting
formulation is discretized in space by finite elements and in time by an
unconditionally stable semi-implicit scheme, and is implemented as a reusable
\code{ImmersedBoundary} physics kernel within the \grins multiphysics framework,
built on the \libmesh finite element library. The method is verified on three
FSI benchmarks---an elliptically displaced thick ring returning to equilibrium,
a radially stretched incompressible ring, and a disk falling under gravity in a
viscous fluid---for which the mixed formulation removes the volumetric failure
observed with the unstabilized formulation and reproduces the analytical
terminal velocity of the falling disk to within~$1\%$. We also report a
pressure-loaded inflating-ring case that exposes the remaining sensitivity of
the method to the solid-to-fluid mesh-size ratio, and discuss it as an open
challenge.
\end{abstract}

\begin{keyword}
immersed boundary method \sep finite element method \sep fluid--structure
interaction \sep fictitious domain \sep distributed Lagrange multiplier \sep
volumetric stabilization \sep incompressible hyperelasticity
\end{keyword}

\end{frontmatter}

\section{Introduction}
\label{sec:intro}

Fluid--structure interaction (FSI) couples two of the most challenging problems
in continuum mechanics and arises across a wide range of engineering
applications, from aircraft and parachute design to cardiovascular and other
biomedical simulations. In many of these applications the structure undergoes
large deformations while remaining fully immersed in the surrounding fluid, and
the ability to predict the coupled response robustly and efficiently remains a
central goal of computational mechanics.

One of the most well established approaches for simulating FSI is the arbitrary
Lagrangian--Eulerian (ALE) framework~\cite{ale}. In the ALE approach the fluid
and the solid share a conforming interface and do not overlap, so that whenever
the solid deforms the fluid mesh must deform with it. For large structural
deformations the fluid elements near the interface become severely distorted,
and the fluid domain must be remeshed. The complexity and cost of this mesh
reconstruction make ALE difficult to automate for problems involving very large
elastic deformations, such as the deployment of a parachute system.

An alternative that avoids remeshing altogether is the immersed boundary method
(IBM). Unlike ALE, in which a single partitioned mesh is used, the IBM
completely immerses the solid in the fluid and employs two overlapping and
independent grids---an Eulerian grid for the fluid and a Lagrangian grid for the
structure---so that no remeshing is required as the structure moves. The central
difficulty is that the fluid and the structure are naturally described in
different mathematical frameworks: the incompressible Navier--Stokes equations
are most conveniently posed in an Eulerian frame in terms of velocity and
pressure, whereas the constitutive response of a hyperelastic solid is most
naturally expressed in a Lagrangian frame in terms of the deformation of
material particles. The IBM links the two frameworks
together~\cite{hyperelastic}, using one grid for the Eulerian fluid variables
and another overlapping grid for the Lagrangian variables of the immersed
structure.

The IBM was introduced by Peskin in the 1970s to simulate blood flow around
heart valves~\cite{peskin} and has since been applied extensively to problems in
biological fluid dynamics. In the original formulation the Navier--Stokes
equations are solved over the entire domain and the presence of the solid is
accounted for through a body force expressed by means of a Dirac delta
distribution~\cite{peskin}. The delta distribution acts as an interpolation
kernel between the Eulerian and Lagrangian frameworks and implicitly performs
the change of variables needed to use both descriptions
simultaneously~\cite{heltaithesis}. The accuracy of the method depends critically
on the construction of a suitable regularized approximation of the delta
distribution, and the original formulation, based on finite differences for both
the fluid and the solid, was moreover restricted to structures that can be
modelled as fibre-like (co-dimension one) bodies~\cite{heltaithesis}.

A significant generalization was provided by the finite element immersed boundary
method (FE-IBM) of Boffi, Gastaldi and co-workers~\cite{boffi2003,hyperelastic}.
By discretizing the coupled problem with finite elements, the delta distribution
term is treated variationally through the test functions rather than
explicitly, which improves the representation of the immersed interface and
extends the method to general (thick, co-dimension zero) hyperelastic
solids~\cite{hyperelastic,boffi2011}. A further modification was proposed
in~\cite{boffilm}, in which a distributed Lagrange multiplier is introduced in
the spirit of the fictitious domain method~\cite{glow1,glow2} to enforce the
kinematic constraint that the velocity of the solid match that of the
surrounding fluid. The multiplier keeps the fluid and solid equations largely
separated and, importantly, a semi-implicit time discretization of the resulting
distributed-Lagrange-multiplier FE-IBM (DLM-FE-IBM) was shown to be
unconditionally stable with respect to the time step~\cite{boffilm}. This
fictitious domain, distributed-Lagrange-multiplier formulation has continued to
be developed actively, including detailed studies of the coupling/interface
matrix and of the non-matching quadrature it requires, parallel
implementations, and a posteriori error analysis, recently surveyed
in~\cite{alshehri2025}.

\begin{figure}[t]
    \centering
    \begin{subfigure}[b]{0.38\textwidth}
        \includegraphics[width=\textwidth]{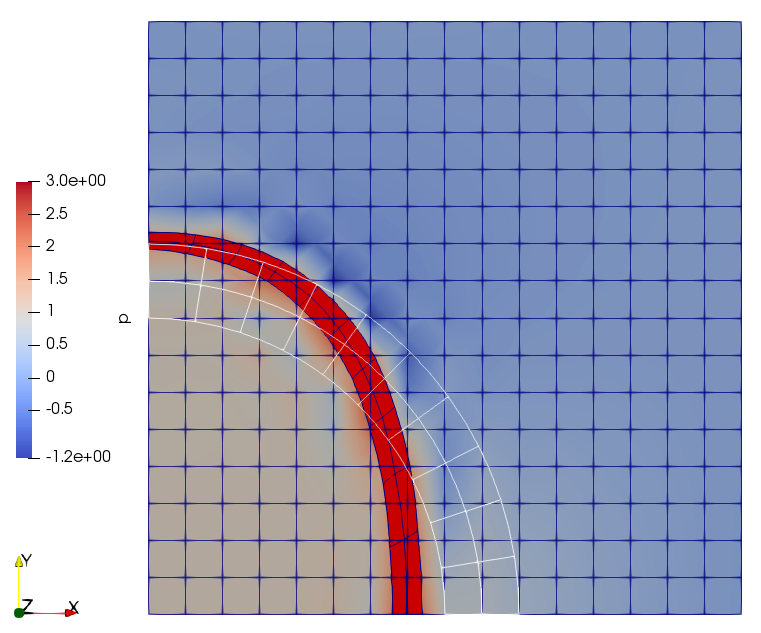}
        \caption{$t = 0.01$ s}
    \end{subfigure}\hspace{0.04\textwidth}%
    \begin{subfigure}[b]{0.38\textwidth}
        \includegraphics[width=\textwidth]{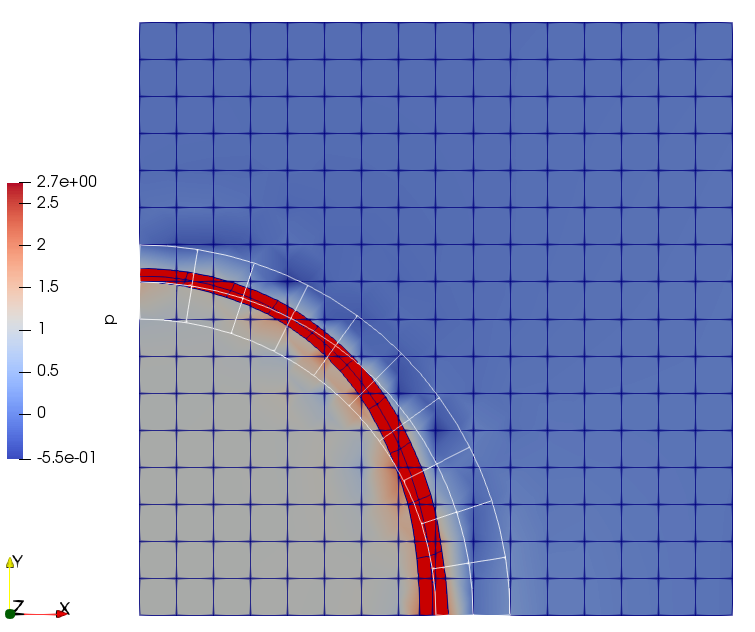}
        \caption{$t = 1$ s}
    \end{subfigure}
    \caption{Volumetric failure of the unstabilized distributed-Lagrange-multiplier
    FE-IBM. An elliptically displaced incompressible hyperelastic ring, which
    should return to its circular equilibrium while preserving its area, instead
    collapses as the system evolves.}
    \label{fig:fail}
\end{figure}

We are interested in the fully variational, finite element form of the IBM
because it relies entirely on finite elements and therefore fits naturally into
the \grins multiphysics framework~\cite{grins}, a \code{C++} library built on the
\libmesh parallel finite element library~\cite{libmesh}. On implementing the
DLM-FE-IBM in \grins and running standard FSI benchmarks, however, we observed
that when the immersed solid is modelled as a \emph{fully incompressible}
hyperelastic material a direct discretization of the deviatoric solid stress
produces severe volumetric instabilities. Figure~\ref{fig:fail} illustrates the
difficulty for a thick ring, initially displaced into an ellipse, that should
relax back to a circular equilibrium: instead of preserving its area the ring
collapses as the system evolves. The underlying issue is well known in
computational solid mechanics---the incompressibility of a hyperelastic
structure is satisfied only approximately after discretization~\cite{boyce1}---
and it is not always visible in the FSI literature because some formulations
employ non-standard fluid spaces or divergence-conforming kernels that
implicitly reinforce incompressibility.

The problem of enforcing incompressibility for immersed hyperelastic structures
has attracted considerable recent attention. Vadala-Roth et
al.~\cite{boyce1} introduced a volumetric stabilization for the immersed
finite-element/finite-difference (IFED) method in which a dilatational energy is
added to the solid region so that the Eulerian incompressibility constraint is
reinforced by additional structural stresses driving the Lagrangian constraint
$J=1$; the stabilization is active only in the solid, is compatible with standard
nodal Lagrangian elements, and yields accuracy comparable to specialized
incompressible-elasticity solvers. More recent work has pursued alternatives
that avoid an explicit volumetric penalty altogether, notably
divergence-free composite B-spline interpolation kernels that preserve the
discrete incompressibility of the transferred velocity
field~\cite{griffithcbs}. Closely related pure-solid stabilized finite element
formulations for nearly and fully incompressible finite elasticity, based on
mixed displacement--pressure spaces or projection stabilization, have also been
developed independently of the immersed setting~\cite{stein,rossi}. A common
thread across this body of work is that treating the solid pressure explicitly,
rather than relying on the fluid pressure alone, is essential for robust volume
conservation of an incompressible immersed structure.

\paragraph{Contributions.} In this work we develop and verify a
volumetrically stabilized \emph{mixed formulation} of the
distributed-Lagrange-multiplier FE-IBM for fully incompressible hyperelastic
solids. Our contributions are the following. First, building on the
nearly-incompressible hyperelasticity framework~\cite{stein} and the
stabilization idea of~\cite{boyce1}, we augment the immersed solid stress with a
volumetric contribution and introduce an additional solid pressure field~$p_s$,
enforced weakly, that acts as the Lagrange multiplier of the Lagrangian
incompressibility constraint. Unlike the IFED stabilization of~\cite{boyce1},
which is formulated within a finite-difference immersed method, our formulation
is posed and discretized entirely within the fully variational DLM-FE-IBM,
making the solid pressure a genuine mixed unknown. Second, we cast the
formulation for a fully incompressible neo-Hookean material, derive its
semi-implicit and spatial finite element discretizations, and describe a reusable
implementation as an \code{ImmersedBoundary} physics kernel in \grins, including
the overlapping fluid--solid map required to assemble the non-matching coupling
terms in parallel. Third, we verify the method on three FSI benchmarks and show
that the mixed formulation removes the volumetric failure of the unstabilized
scheme, and we report an inflating-ring test case that exposes a remaining
sensitivity to the solid-to-fluid mesh-size ratio, which we discuss as an open
challenge.

The remainder of the paper is organized as follows.
Section~\ref{sec:formulation} summarizes the governing equations and the existing
IBM formulations, derives the proposed mixed formulation, and presents its time
and space discretizations. Section~\ref{sec:implementation} describes the
implementation within the \grins framework. Section~\ref{sec:results} presents
the numerical experiments, and Section~\ref{sec:conclusions} concludes and
outlines directions for future work.

\section{Governing equations and the mixed immersed boundary formulation}
\label{sec:formulation}

In this section we first fix notation and recall the continuum-mechanical
ingredients needed for the immersed boundary method
(Section~\ref{sec:prelim}). We then summarize, in a common variational setting,
the existing FE-IBM formulations---the hyperelastic and variational
formulations~\cite{hyperelastic} and the distributed-Lagrange-multiplier
fictitious domain formulation~\cite{boffilm}
(Section~\ref{sec:existing})---so that the paper is self-contained. Our
contribution begins in Section~\ref{sec:mixed}, where we derive the proposed
mixed formulation for a fully incompressible neo-Hookean solid. Its temporal and
spatial discretizations are given in Sections~\ref{sec:time}
and~\ref{sec:space}.

\subsection{Preliminaries and notation}
\label{sec:prelim}

Let $B \subset \mathbb{R}^d$, $d\in\{2,3\}$, denote the reference (undeformed)
configuration of the immersed solid and let $\Omega \subset \mathbb{R}^d$ denote
the fixed control volume occupied by the fluid--solid system. Points in the
reference configuration are denoted by $\bv{X}$ and points in the deformed
configuration by $\bv{x}$. The motion of the solid is described by the
deformation map $\bv{\varphi}: B \times [0,T] \to \Omega$, which maps each
material point $\bv{X}\in B$ to its current position
$\bv{x} = \bv{\varphi}(\bv{X},t)$, and $B_t = \bv{\varphi}(B,t)$ denotes the
current configuration of the solid. Eulerian fields are functions of~$\bv{x}$
and Lagrangian fields functions of~$\bv{X}$. The deformation gradient and
displacement are
\begin{equation}\label{eqn:F}
\bv{F}(\bv{X},t) = \nabla_{_{\bv{X}}}\bv{\varphi}(\bv{X},t),
\qquad
\bv{u}(\bv{X},t) = \bv{\varphi}(\bv{X},t) - \bv{X},
\end{equation}
with $J = \det\bv{F}$, and the velocity satisfies
$\bv{v}(\bv{x},t) = \partial\bv{\varphi}/\partial t$ for
$\bv{x} = \bv{\varphi}(\bv{X},t)$. Following~\cite{hyperelastic} the Cauchy
stress of the immersed system is additively split as
\begin{equation}\label{eqn:Cauchy}
\bv{\sigma} =
\begin{cases}
\bv{\sigma_f}, & \text{in } \Omega\backslash B_t,\\[2pt]
\bv{\sigma_f} + \bv{\sigma_s}, & \text{in } B_t,
\end{cases}
\qquad
\rho =
\begin{cases}
\rho_f, & \text{in } \Omega\backslash B_t,\\[2pt]
\rho_s, & \text{in } B_t,
\end{cases}
\end{equation}
which expresses that the solid is immersed in, and behaves as an addition to,
the surrounding fluid; the mass densities $\rho_f,\rho_s$ are piecewise constant
with $\rho_s \ge \rho_f$, and we write $\delta\rho = \rho_s - \rho_f$. The fluid
is an incompressible Newtonian fluid, so that
\begin{equation}\label{eqn:sigmaf}
\bv{\sigma_f} = -p\bv{I} + \mu\bigl(\grad\bv{v} + (\grad\bv{v})^T\bigr),
\qquad
\diverg{\bv{v}} = 0 ,
\end{equation}
with dynamic viscosity~$\mu$ and pressure~$p$. The solid stress is expressed
through the first Piola--Kirchhoff stress $\bv{P} = J\bv{\sigma_s}\bv{F}^{-T}$.
For a hyperelastic material with strain energy density $W$,
\begin{equation}\label{eqn:PW}
\bv{P}(\bv{F}) = \pdv{W}{\bv{F}},
\qquad
\bv{S} = \bv{F}^{-1}\bv{P} = 2\,\pdv{W}{\bv{C}},
\qquad
\bv{C} = \bv{F}^T\bv{F},
\end{equation}
where $\bv{S}$ is the second Piola--Kirchhoff stress and $\bv{C}$ the right
Cauchy--Green tensor. We use the standard velocity and pressure spaces
$\bv{V} = (H^1_0(\Omega))^d$ and $Q = L^2_0(\Omega)$, denote by $(\cdot,\cdot)$
the $L^2$ inner product and by $\langle\cdot,\cdot\rangle$ the duality pairing,
and introduce the bilinear and trilinear forms~\cite{gunz}
\begin{equation}\label{eqn:forms}
a(\bv{p},\bv{q}) = \mu\bigl(\nabla_{sym}\bv{p},\nabla_{sym}\bv{q}\bigr),
\qquad
b(\bv{p},\bv{q},\bv{r}) = \tfrac{\rho_f}{2}
\bigl[(\bv{p}\cdot\grad\bv{q},\bv{r}) - (\bv{p}\cdot\grad\bv{r},\bv{q})\bigr],
\end{equation}
with $\nabla_{sym}\bv{p} = \grad\bv{p} + (\grad\bv{p})^T$.

\subsection{The finite element immersed boundary method}
\label{sec:existing}

Starting from the principle of virtual work with the stress and density splits
\eqref{eqn:Cauchy}, and expressing the transmission of the solid elastic
response to the fluid through the delta distribution, Boffi et
al.~\cite{hyperelastic,boffi2011} derived the hyperelastic FE-IBM. Its
variational form avoids any explicit construction of the delta distribution:
the elastic forces enter through the test functions evaluated along the
deformation map. Writing $\bv{\phi}$ for a fluid test function, the variational
FE-IBM reads: for all $t\in[0,T]$, find
$(\bv{v}(t),p(t)) \in \bv{V}\times Q$ and
$\bv{\varphi}(t): B \to \Omega$ such that
\begin{subequations}\label{eqn:varibm}
\begin{align}
&\rho_f\tfrac{d}{dt}(\bv{v},\bv{\phi}) + b(\bv{v},\bv{v},\bv{\phi})
   + a(\bv{v},\bv{\phi}) - (\diverg\bv{\phi},p)
   = \langle\bv{d},\bv{\phi}\rangle + \langle\bv{f},\bv{\phi}\rangle,
   && \forall\,\bv{\phi}\in\bv{V},\\
&(\diverg\bv{v},q) = 0, && \forall\,q\in Q,\\
&\langle\bv{f},\bv{\phi}\rangle
   = -\!\int_B \bv{P}(\bv{F}) : \nabla_{_{\bv{X}}}\bv{\phi}(\bv{\varphi})\,d\bv{X},
   \quad
   \langle\bv{d},\bv{\phi}\rangle
   = -\delta\rho\!\int_B \tfrac{\partial^2\bv{\varphi}}{\partial t^2}
       \cdot\bv{\phi}(\bv{\varphi})\,d\bv{X}, &&\\
&\pdv{\bv{\varphi}}{t}(\bv{X},t) = \bv{v}(\bv{\varphi}(\bv{X},t),t),
   && \forall\,\bv{X}\in B,
\end{align}
\end{subequations}
together with the initial conditions $\bv{v}(\cdot,0)=\bv{v}_0$ and
$\bv{\varphi}(\cdot,0)=\bv{\varphi}_0$. The term $\bv{f}$ carries the internal
elastic force of the solid and $\bv{d}$ accounts for the density difference
between solid and fluid.

To decouple the fluid and solid equations and obtain a scheme that is
unconditionally stable in time, Boffi et al.~\cite{boffilm} enforced the
kinematic constraint (\ref{eqn:varibm}d) weakly through a distributed Lagrange
multiplier $\bv{\lambda}$, in the spirit of the fictitious domain
method~\cite{glow1,glow2}. Setting $\bv{\Lambda} = (H^1(B))^d$ and defining the
coupling form
\begin{equation}\label{eqn:cform}
\bv{c}(\bv{\lambda},\bv{Z})
= \int_B \bigl(\nabla_{_{\bv{X}}}\bv{\lambda} : \nabla_{_{\bv{X}}}\bv{Z}
   + \bv{\lambda}\cdot\bv{Z}\bigr)\,d\bv{X},
\qquad \bv{\lambda},\bv{Z}\in\bv{\Lambda},
\end{equation}
the distributed-Lagrange-multiplier FE-IBM (DLM-FE-IBM) reads: find
$(\bv{v},p)\in\bv{V}\times Q$, $\bv{\varphi}\in (H^1(B))^d$ and
$\bv{\lambda}\in\bv{\Lambda}$ such that
\begin{subequations}\label{eqn:fdibm}
\begin{align}
&\rho_f\tfrac{d}{dt}(\bv{v},\bv{\phi_f}) + b(\bv{v},\bv{v},\bv{\phi_f})
   + a(\bv{v},\bv{\phi_f}) - (\diverg\bv{\phi_f},p)
   + \bv{c}(\bv{\lambda},\bv{\phi_f}(\bv{\varphi})) = 0,
   && \forall\,\bv{\phi_f}\in\bv{V},\\
&(\diverg\bv{v},q) = 0, && \forall\,q\in Q,\\
&\delta\rho\Bigl(\tfrac{\partial^2\bv{\varphi}}{\partial t^2},\bv{\phi_s}\Bigr)_B
   + \bigl(\bv{P}(\bv{F}),\nabla_{_{\bv{X}}}\bv{\phi_s}\bigr)_B
   - \bv{c}(\bv{\lambda},\bv{\phi_s}) = 0,
   && \forall\,\bv{\phi_s}\in (H^1(B))^d,\\
&\bv{c}\Bigl(\bv{\phi_\lambda},\,\bv{v}(\bv{\varphi}(\bv{X},t),t)
   - \pdv{\bv{\varphi}}{t}(\bv{X},t)\Bigr) = 0,
   && \forall\,\bv{\phi_\lambda}\in\bv{\Lambda},
\end{align}
\end{subequations}
with $(\cdot,\cdot)_B$ the $L^2(B)^d$ inner product. The multiplier
$\bv{\lambda}$ enforces the velocity-matching constraint and mediates the
transfer of the solid elastic response to the fluid, keeping the two subproblems
largely separated.

\subsection{The proposed mixed formulation}
\label{sec:mixed}

Formulation~\eqref{eqn:fdibm} enforces incompressibility of the coupled system
only through the Eulerian constraint $\diverg\bv{v}=0$ acting on the fluid
velocity. When the immersed solid is a \emph{fully incompressible} hyperelastic
material this is not sufficient: after discretization the Lagrangian
incompressibility of the structure, $J = \det\bv{F} = 1$, is satisfied only
approximately, and the deviatoric elastic stress alone cannot prevent the
spurious volume changes and locking illustrated in Fig.~\ref{fig:fail}. This is
the immersed counterpart of the classical volumetric locking of incompressible
elasticity~\cite{boyce1,stein}.

To restore volumetric stability we introduce a solid pressure field that acts as
a Lagrange multiplier for the constraint $J=1$ within the solid. In continuum
mechanics the Cauchy stress admits the deviatoric/volumetric
decomposition~\cite{cmbook}
\begin{equation}\label{eqn:dev}
\bv{\sigma} = \operatorname{dev}[\bv{\sigma}] + p\,\bv{I},
\end{equation}
in which, for incompressible motion, $p$ is precisely the Lagrange multiplier
enforcing incompressibility. In the immersed setting the fluid-like stress
$\bv{\sigma_f}$ already carries such a term through the fluid pressure, but the
solid stress $\bv{\sigma_s}$ is not, in general, deviatoric. We therefore
augment the stress split~\eqref{eqn:Cauchy} with an additional volumetric
contribution acting only in the solid,
\begin{equation}\label{eqn:newcauchy}
\bv{\sigma} = \operatorname{dev}[\bv{\sigma_f}] + p\,\bv{I} +
\begin{cases}
0, & \text{in } \Omega\backslash B_t,\\[2pt]
\bv{\sigma_s} + p_s\bv{I}, & \text{in } B_t,
\end{cases}
\end{equation}
where $p_s$ is a solid pressure that reinforces incompressibility in the solid
region. Consistently, and following the mixed treatment of nearly incompressible
hyperelasticity~\cite{stein}, we split the strain energy into volumetric and
isochoric parts,
\begin{equation}\label{eqn:Wdecomp}
W(\bv{C}) = \kappa\,U(J) + \widetilde{W}(\bv{C}),
\end{equation}
where $\kappa$ is a bulk modulus independent of the deformation, $U(J)$ is a
purely volumetric energy depending only on $J$, and $\widetilde{W}$ is the
isochoric part. The solid pressure is then defined by
\begin{equation}\label{eqn:ps}
p_s = \kappa\,U'(J),
\end{equation}
so that the volumetric stabilization~\eqref{eqn:newcauchy} exactly parallels the
mixed formulation of a nearly incompressible hyperelastic material when the
solid is not immersed~\cite{stein}.

Using~\eqref{eqn:PW} and~\eqref{eqn:Wdecomp}, the second Piola--Kirchhoff stress
becomes
\begin{equation}\label{eqn:Snew}
\bv{S} = 2\kappa\,U'(J)\,\pdv{J}{\bv{C}} + 2\,\pdv{\widetilde{W}}{\bv{C}} .
\end{equation}
For a neo-Hookean isochoric energy
$\widetilde{W}(\bv{C}) = \tfrac{\mu_s}{2}(J^{-2/3}I_1 - 3)$, with
$I_1 = \operatorname{tr}\bv{C}$, one has
\begin{equation}\label{eqn:dWdC}
\pdv{\widetilde{W}}{\bv{C}}
= \frac{\mu_s}{2}\,J^{-2/3}\Bigl(\bv{I} - \tfrac{I_1}{3}\bv{C}^{-1}\Bigr).
\end{equation}
Substituting~\eqref{eqn:ps},~\eqref{eqn:dWdC} and the identity
$\partial J/\partial\bv{C} = \tfrac{1}{2}J\bv{C}^{-1}$ into~\eqref{eqn:Snew}, and
mapping to the first Piola--Kirchhoff stress, gives the stabilized solid stress
\begin{equation}\label{eqn:Pnew}
\bv{P}(\bv{F}) = p_s\,J\,\bv{F}^{-T}
   + \mu_s\,J^{-2/3}\Bigl(\bv{F} - \tfrac{I_1}{3}\bv{F}^{-T}\Bigr).
\end{equation}
The first term is the volumetric stabilization; the second is the standard
isochoric neo-Hookean response.

The definition~\eqref{eqn:ps} is imposed weakly. For a nearly incompressible
material this reads
$\int_B\bigl[U'(J) - p_s/\kappa\bigr]\,q_s\,d\bv{X} = 0$ for all
$q_s\in L^2(B)$. The fully incompressible limit is recovered as
$\kappa\to\infty$, i.e.\ $1/\kappa\to 0$, giving the constraint
\begin{equation}\label{eqn:Uconstraint}
\int_B U'(J)\,q_s\,d\bv{X} = 0
\qquad \forall\,q_s\in L^2(B),
\end{equation}
in which $p_s$ has become a genuine Lagrange multiplier. Following~\cite{boyce1}
we use the volumetric energy
\begin{equation}\label{eqn:U}
U(J) = \tfrac{1}{2}(\ln J)^2,
\qquad U'(J) = \frac{\ln J}{J},
\end{equation}
which penalizes deviations of $J$ from unity and vanishes only at $J=1$.

Collecting~\eqref{eqn:newcauchy}--\eqref{eqn:U} and adding the solid pressure as
a mixed unknown to the DLM-FE-IBM~\eqref{eqn:fdibm} yields the proposed
\emph{mixed formulation} for a fully incompressible neo-Hookean solid immersed
in an incompressible Newtonian fluid: find
$(\bv{v},p)\in\bv{V}\times Q$,
$(\bv{\varphi},p_s)\in (H^1(B))^d \times L^2(B)$ and
$\bv{\lambda}\in\bv{\Lambda}$ such that
\begin{subequations}\label{eqn:mixedibm}
\begin{align}
&\rho_f\tfrac{d}{dt}(\bv{v},\bv{\phi_f}) + b(\bv{v},\bv{v},\bv{\phi_f})
   + a(\bv{v},\bv{\phi_f}) - (\diverg\bv{\phi_f},p)
   + \bv{c}(\bv{\lambda},\bv{\phi_f}(\bv{\varphi})) = 0,
   && \forall\,\bv{\phi_f}\in\bv{V},\\
&(\diverg\bv{v},q) = 0, && \forall\,q\in Q,\\
&\delta\rho\Bigl(\tfrac{\partial^2\bv{\varphi}}{\partial t^2},\bv{\phi_s}\Bigr)_B
   + \bigl(\bv{P}(\bv{F}),\nabla_{_{\bv{X}}}\bv{\phi_s}\bigr)_B
   - \bv{c}(\bv{\lambda},\bv{\phi_s}) = 0,
   && \forall\,\bv{\phi_s}\in (H^1(B))^d,\\
&\bigl(U'(J(\bv{F})),q_s\bigr)_B = 0, && \forall\,q_s\in L^2(B),\\
&\bv{c}\Bigl(\bv{\phi_\lambda},\,\bv{v}(\bv{\varphi}(\bv{X},t),t)
   - \pdv{\bv{\varphi}}{t}(\bv{X},t)\Bigr) = 0,
   && \forall\,\bv{\phi_\lambda}\in\bv{\Lambda},
\end{align}
\end{subequations}
with $\bv{P}(\bv{F})$ given by~\eqref{eqn:Pnew}, $U$ by~\eqref{eqn:U}, and the
initial conditions $\bv{v}(\cdot,0)=\bv{v}_0$, $\bv{\varphi}(\cdot,0)=
\bv{\varphi}_0$. Compared with~\eqref{eqn:fdibm}, the mixed
formulation~\eqref{eqn:mixedibm} adds the solid pressure unknown $p_s$ and the
weak constraint (\ref{eqn:mixedibm}d), and replaces the deviatoric solid stress
by the stabilized stress~\eqref{eqn:Pnew}. In contrast to the finite-difference
IFED stabilization of~\cite{boyce1}, here $p_s$ is a mixed finite element field
living on the solid mesh, so that the volumetric constraint is enforced within
the same fully variational setting as the fluid pressure and the distributed
multiplier.

\subsection{Time discretization}
\label{sec:time}

We partition $[0,T]$ into $N$ steps of size $\Delta t$, with $t_n = n\Delta t$
and superscripts denoting time levels. A fully implicit (backward-Euler)
discretization of~\eqref{eqn:mixedibm} requires, at each Newton iteration, the
evaluation of the fluid shape functions along the updated deformation
$\bv{\varphi}^{n+1}$ through the coupling term
$\bv{c}(\bv{\lambda}^{n+1},\bv{\phi_f}(\bv{\varphi}^{n+1}))$. Because these
evaluations must be recomputed for every linearization, the fully implicit scheme
is computationally expensive.

Following~\cite{boffilm,heltaithesis} we adopt instead a semi-implicit
(forward-Euler/backward-Euler) scheme in which the solid position entering the
coupling term is taken at the previous time level, so that the fluid shape
functions are evaluated at the known configuration $\bv{\varphi}^{n}$. The
elastic response of the solid (the first Piola--Kirchhoff stress) is still
evaluated at the current configuration, but the corresponding body force is
applied to the fluid at the location of the solid at the previous
step~\cite{heltaiibm}. This semi-implicit discretization of the DLM-FE-IBM is
unconditionally stable with respect to the time step~\cite{boffilm}. The scheme
reads: for $n=1,\dots,N$, find $(\bv{v}^{n+1},p^{n+1})\in\bv{V}\times Q$,
$(\bv{\varphi}^{n+1},p_s^{n+1})\in (H^1(B))^d\times L^2(B)$ and
$\bv{\lambda}^{n+1}\in\bv{\Lambda}$ such that
\begin{subequations}\label{eqn:semiimplicit}
\begin{align}
&\rho_f\Bigl(\tfrac{\bv{v}^{n+1}-\bv{v}^{n}}{\Delta t},\bv{\phi_f}\Bigr)
   + b(\bv{v}^{n+1},\bv{v}^{n+1},\bv{\phi_f}) + a(\bv{v}^{n+1},\bv{\phi_f})
   \nonumber\\
&\qquad - (\diverg\bv{\phi_f},p^{n+1})
   + \bv{c}(\bv{\lambda}^{n+1},\bv{\phi_f}(\bv{\varphi}^{n})) = 0,
   && \forall\,\bv{\phi_f}\in\bv{V},\\
&(\diverg\bv{v}^{n+1},q) = 0, && \forall\,q\in Q,\\
&\delta\rho\Bigl(\tfrac{\bv{\varphi}^{n+1}-2\bv{\varphi}^{n}+\bv{\varphi}^{n-1}}
   {\Delta t^2},\bv{\phi_s}\Bigr)_B
   + \bigl(\bv{P}(\bv{F}^{n+1}),\nabla_{_{\bv{X}}}\bv{\phi_s}\bigr)_B
   - \bv{c}(\bv{\lambda}^{n+1},\bv{\phi_s}) = 0,
   && \forall\,\bv{\phi_s}\in (H^1(B))^d,\\
&\bigl(U'(J(\bv{F}^{n+1})),q_s\bigr)_B = 0, && \forall\,q_s\in L^2(B),\\
&\bv{c}\Bigl(\bv{\phi_\lambda},\,\bv{v}^{n+1}(\bv{\varphi}^{n})
   - \tfrac{\bv{\varphi}^{n+1}-\bv{\varphi}^{n}}{\Delta t}\Bigr) = 0,
   && \forall\,\bv{\phi_\lambda}\in\bv{\Lambda}.
\end{align}
\end{subequations}

\subsection{Spatial discretization}
\label{sec:space}

We introduce shape-regular meshes of $\Omega$ and $B$ with characteristic sizes
$h_x$ and $h_s$, respectively. Let $V_h\subseteq\bv{V}$ and $Q_h\subseteq Q$ be
finite element spaces satisfying the discrete inf--sup condition for the
Navier--Stokes problem~\cite{gunz}, and let $U_h\subseteq (H^1(B))^d$,
$S_h\subseteq L^2(B)$ and $\Lambda_h\subseteq\bv{\Lambda}$ be finite element
spaces on the solid mesh. The fully discrete mixed formulation is obtained by
restricting the semi-implicit scheme~\eqref{eqn:semiimplicit} to these spaces:
for $n=1,\dots,N$, find $(\bv{v}^{n+1}_h,p^{n+1}_h)\in V_h\times Q_h$,
$(\bv{\varphi}^{n+1}_h,p^{n+1}_{s,h})\in U_h\times S_h$ and
$\bv{\lambda}^{n+1}_h\in\Lambda_h$ such that
\begin{subequations}\label{eqn:space}
\begin{align}
&\rho_f\Bigl(\tfrac{\bv{v}^{n+1}_h-\bv{v}^{n}_h}{\Delta t},\bv{\phi_f}\Bigr)
   + b(\bv{v}^{n+1}_h,\bv{v}^{n+1}_h,\bv{\phi_f}) + a(\bv{v}^{n+1}_h,\bv{\phi_f})
   \nonumber\\
&\qquad - (\diverg\bv{\phi_f},p^{n+1}_h)
   + \bv{c}(\bv{\lambda}^{n+1}_h,\bv{\phi_f}(\bv{\varphi}^{n}_h)) = 0,
   && \forall\,\bv{\phi_f}\in V_h,\\
&(\diverg\bv{v}^{n+1}_h,q) = 0, && \forall\,q\in Q_h,\\
&\delta\rho\Bigl(\tfrac{\bv{\varphi}^{n+1}_h-2\bv{\varphi}^{n}_h+\bv{\varphi}^{n-1}_h}
   {\Delta t^2},\bv{\phi_s}\Bigr)_B
   + \bigl(\bv{P}(\bv{F}^{n+1}_h),\nabla_{_{\bv{X}}}\bv{\phi_s}\bigr)_B
   - \bv{c}(\bv{\lambda}^{n+1}_h,\bv{\phi_s}) = 0,
   && \forall\,\bv{\phi_s}\in U_h,\\
&\bigl(U'(J(\bv{F}^{n+1}_h)),q_s\bigr)_B = 0, && \forall\,q_s\in S_h,\\
&\bv{c}\Bigl(\bv{\phi_\lambda},\,\bv{v}^{n+1}_h(\bv{\varphi}^{n}_h)
   - \tfrac{\bv{\varphi}^{n+1}_h-\bv{\varphi}^{n}_h}{\Delta t}\Bigr) = 0,
   && \forall\,\bv{\phi_\lambda}\in\Lambda_h.
\end{align}
\end{subequations}
In all computations reported below the fluid velocity/pressure pair
$(V_h,Q_h)$ is an inf--sup stable Taylor--Hood element, and the solid
displacement, solid pressure and multiplier are discretized with continuous
piecewise-polynomial elements on the solid mesh. The coupling
form~\eqref{eqn:cform} involves integrals of fluid shape functions evaluated at
solid quadrature points that lie in different background elements; the accurate
evaluation of these non-matching integrals, and the associated quadrature, is
discussed in Section~\ref{sec:implementation}.

\section{Implementation in the \grins framework}
\label{sec:implementation}

The fully variational structure of the mixed formulation~\eqref{eqn:space} makes
it a natural fit for a general finite element framework. We implement it in
\grins~\cite{grins}, a \code{C++} multiphysics library built on the
\libmesh~\cite{libmesh} parallel adaptive finite element library, and in
particular on the \code{libMesh::FEMSystem} application framework. \grins
provides an extensible mechanism for assembling element-level weak-form
residuals and Jacobians for coupled systems of partial differential equations
and for reusing them across problems: physics, materials, variables, boundary
and initial conditions are all selected at run time through an input file, and a
given physics kernel can be restricted to a subdomain of the mesh. This
subdomain-restricted design is what allows an incompressible-flow physics to be
active everywhere while the immersed-solid physics is active only on the solid
region.

\subsection{The \code{ImmersedBoundary} physics kernel}

We implement the mixed formulation as a new \code{ImmersedBoundary} physics
kernel, enabled on the solid subdomain, together with an existing
incompressible-flow physics (\code{Stokes} or \code{IncompressibleNavierStokes})
enabled on the fluid subdomain. The fluid residual and Jacobian contributions
are reused unchanged from the existing incompressible-flow kernels; the
\code{ImmersedBoundary} kernel adds the solid momentum
residual~(\ref{eqn:space}c), the stabilized stress~\eqref{eqn:Pnew}, the weak
volumetric constraint~(\ref{eqn:space}d) for the solid pressure~$p_s$, and the
coupling terms~(\ref{eqn:space}a) and~(\ref{eqn:space}e) involving the
distributed multiplier~$\bv{\lambda}$. Fluid and solid materials---and hence the
densities and moduli $\rho_f,\rho_s,\mu,\mu_s$---are specified independently, so
the immersed solid need not share the density or the viscous response of the
surrounding fluid.

\subsection{The overlapping fluid--solid map}

The distinctive computational element of the immersed formulation is the coupling
between the fixed background fluid mesh and the moving foreground solid mesh,
which is required by the coupling form~\eqref{eqn:cform} and the elastic force
transfer. This is handled by an \code{OverlappingFluidSolidMap} object that
records, for each solid element, the background fluid elements it overlaps
together with the indices of the solid quadrature points contained in each,
and the inverse map from fluid to solid elements. Using this map, the kernel
loops over the solid quadrature points falling within each overlapped fluid
element and assembles the non-matching coupling integrals of~\eqref{eqn:space}
directly, without repeated geometric point-location searches. Because the solid
moves, the map is rebuilt at every time step for the semi-implicit scheme (and
would be rebuilt at every Newton step for a fully implicit scheme).

On a distributed-memory parallel mesh the fluid and solid elements that overlap
may reside on different processors. A \code{parallel\_sync} operation packs, for
each solid element, the pairs of solid and fluid element indices and communicates
them to the processor owning the corresponding fluid element, so that the
association is available locally without redundant searching. To obtain the
correct global sparsity pattern for the coupled Jacobian, an
\code{ImmersedBoundaryCouplingFunctor}, derived from
\code{libMesh::GhostingFunctor}, augments the sparsity pattern with the coupling
blocks between overlapping fluid and solid degrees of freedom and provides the
algebraic ghosting needed on the distributed mesh; the sparsity pattern is
rebuilt together with the overlapping map.

Because the coupling integrals evaluate piecewise-polynomial fluid shape
functions at solid quadrature points that generally lie in the interior of
background elements, the integrands are only piecewise smooth over a solid
element and the standard quadrature is inexact. As is well
documented for non-matching immersed
discretizations, this quadrature error must be controlled to recover the
expected accuracy. In practice we supply additional quadrature during assembly
and study its effect in Section~\ref{sec:results}; consistent with the
non-matching integration analyses of the fictitious-domain FE-IBM, we find that a
moderate amount of extra quadrature is sufficient to obtain converged results.

\section{Numerical experiments}
\label{sec:results}

We now verify the mixed formulation~\eqref{eqn:space} on a set of
fluid--structure interaction benchmarks. Sections~\ref{sec:ellipse}
to~\ref{sec:falling} present three cases for which the formulation performs as
expected: an elliptically displaced thick ring relaxing to equilibrium, a
radially stretched incompressible ring, and a disk falling under gravity in a
viscous fluid. Section~\ref{sec:inflate} then reports an inflating-ring test that
exposes a remaining sensitivity of the method and which we regard as an open
challenge. Throughout, the immersed solid is a fully incompressible neo-Hookean
material and the fluid is an incompressible Newtonian fluid; unless stated
otherwise the assembly uses extra quadrature of order~$8$
(see Sections~\ref{sec:stretch} and~\ref{sec:falling} for the sensitivity to this
choice).

\subsection{Elliptically displaced thick ring}
\label{sec:ellipse}

We first consider the classical benchmark of a thick ring displaced into an
ellipse that should relax back to a circular equilibrium while preserving its
area, inspired by the example in~\cite{boffilm}. A circular reference ring is
centred at the origin of a fluid domain $[-L,L]^2$ with $L = 0.5$~mm; the ring
has outer and inner radii $0.3125$~mm and $0.25$~mm. Exploiting symmetry, the
computational domain is reduced to one quarter of the full domain. The ring is
displaced by stretching it $25\%$ vertically and contracting it by the same
amount horizontally, so as to preserve its initial area~\cite{hyperelastic}. The
elastic constant is $\mu_s = 10^4\ \mathrm{N/mm^2}$, and the fluid density and
viscosity are $\rho_f = 1\ \mathrm{kg/mm^3}$ and
$\mu_f = 1\ \mathrm{N\cdot s/mm^2}$. The displaced ring is immersed in a fluid at
rest and the system is allowed to evolve for~$1$~s.

\begin{figure}[t]
    \centering
    \begin{subfigure}[b]{0.32\textwidth}
        \includegraphics[width=\textwidth]{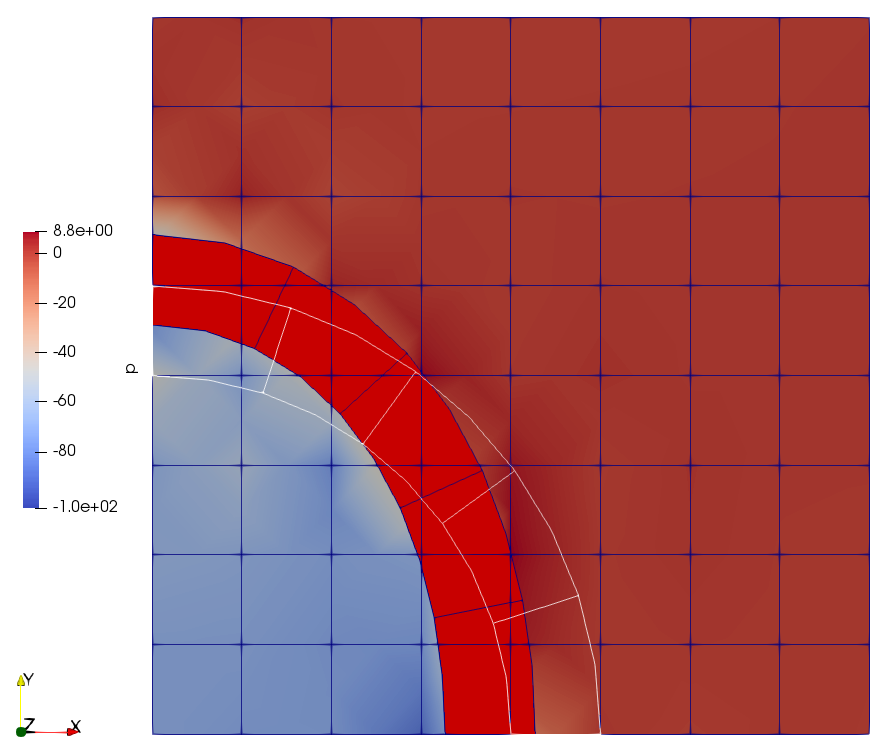}
        \caption{Level 1: $t = 0.01$ s}
    \end{subfigure}%
    \begin{subfigure}[b]{0.32\textwidth}
        \includegraphics[width=\textwidth]{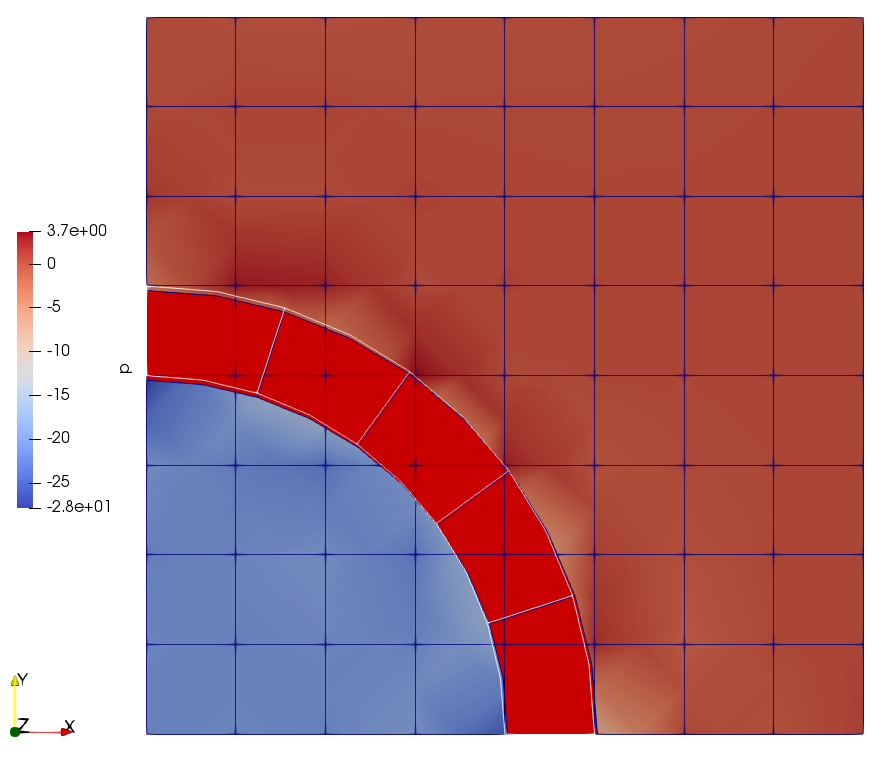}
        \caption{Level 1: $t = 0.13$ s}
    \end{subfigure}%
    \begin{subfigure}[b]{0.32\textwidth}
        \includegraphics[width=\textwidth]{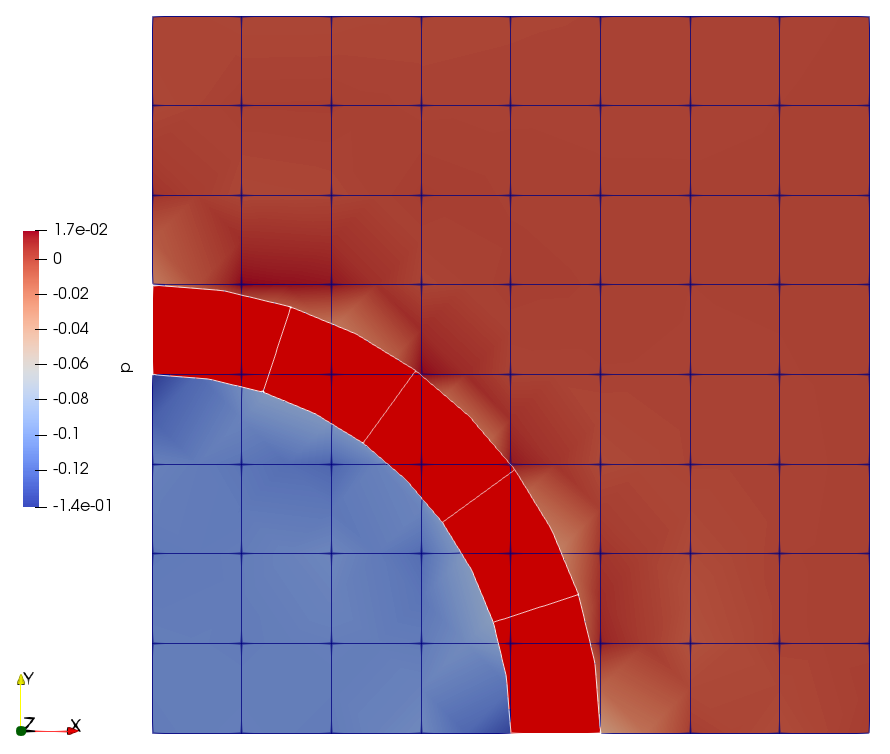}
        \caption{Level 1: $t = 0.3$ s}
    \end{subfigure}

    \begin{subfigure}[b]{0.32\textwidth}
        \includegraphics[width=\textwidth]{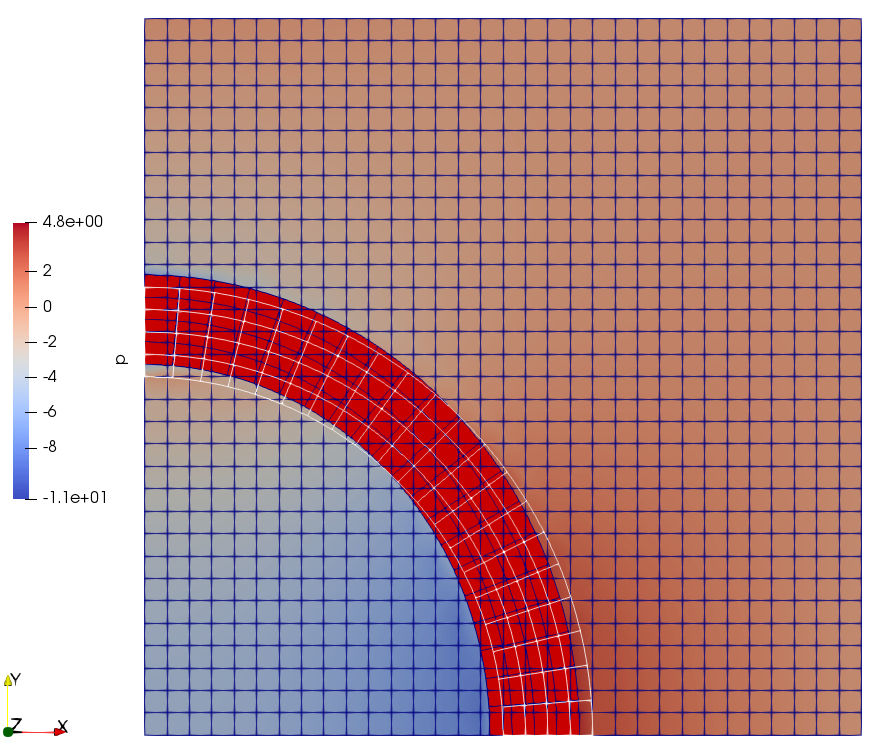}
        \caption{Level 3: $t = 0.01$ s}
    \end{subfigure}%
    \begin{subfigure}[b]{0.32\textwidth}
        \includegraphics[width=\textwidth]{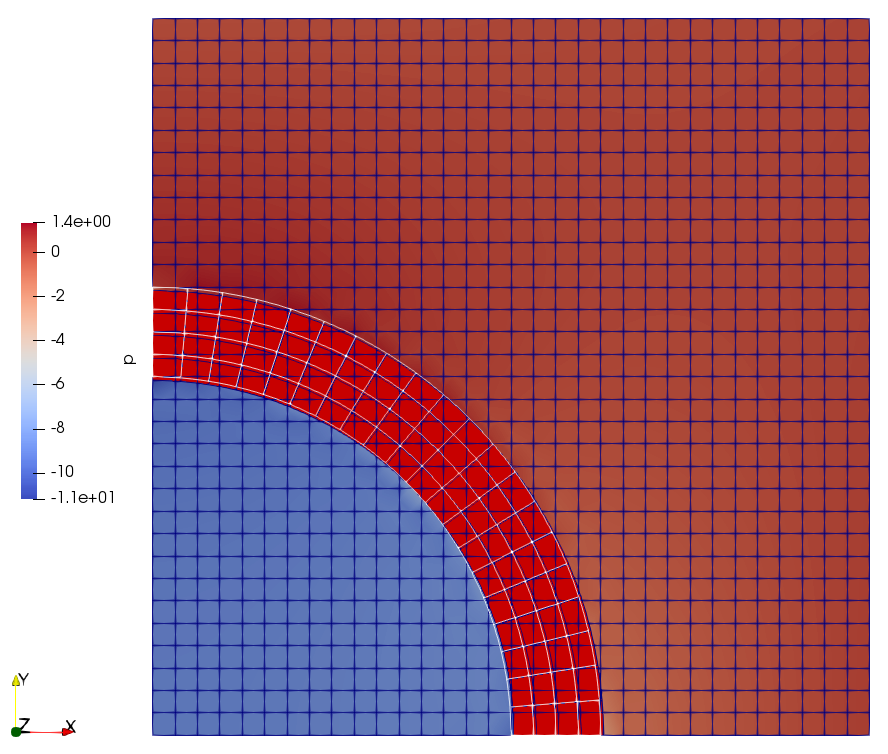}
        \caption{Level 3: $t = 0.07$ s}
    \end{subfigure}%
    \begin{subfigure}[b]{0.32\textwidth}
        \includegraphics[width=\textwidth]{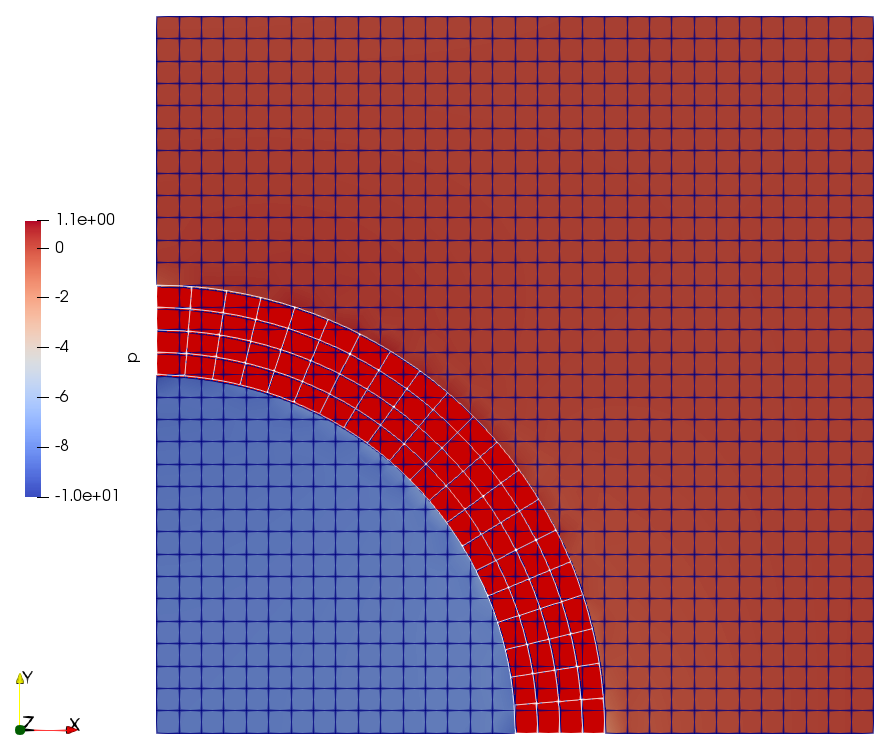}
        \caption{Level 3: $t = 0.11$ s}
    \end{subfigure}
    \caption{Evolution of the elliptically displaced thick ring from its initial
    ellipse to the circular equilibrium, for the coarsest (Level~1, top) and
    finest (Level~3, bottom) meshes of Table~\ref{table:ellipse}. The mixed
    formulation preserves the area and recovers the equilibrium, in contrast to
    the unstabilized result of Fig.~\ref{fig:fail}.
    Here $\mu_f = 1\ \mathrm{N\cdot s/mm^2}$ and $\mu_s = 10^4\ \mathrm{N/mm^2}$.}
    \label{fig:ellipse}
\end{figure}

\begin{figure}[t]
    \centering
    \begin{subfigure}[b]{0.32\textwidth}
        \includegraphics[width=\textwidth]{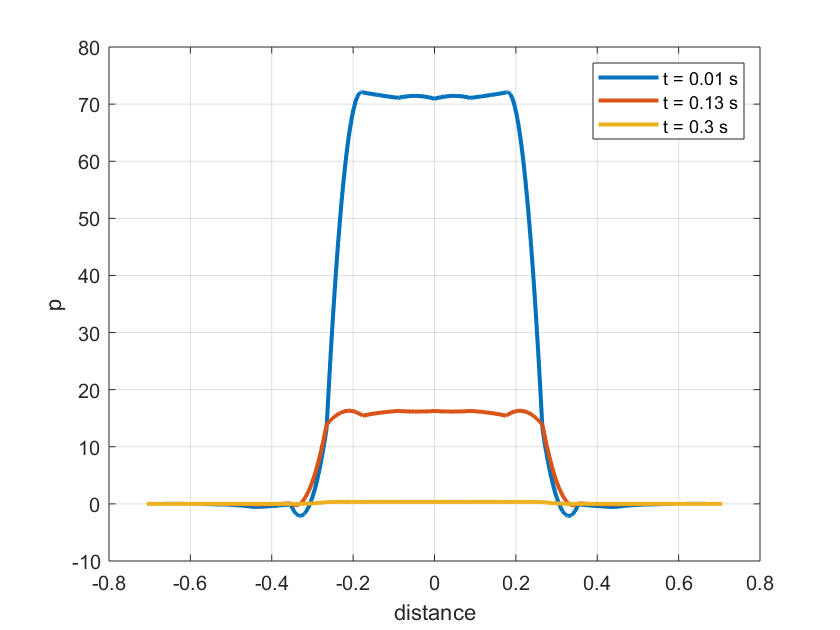}
        \caption{Level 1}
    \end{subfigure}%
    \begin{subfigure}[b]{0.32\textwidth}
        \includegraphics[width=\textwidth]{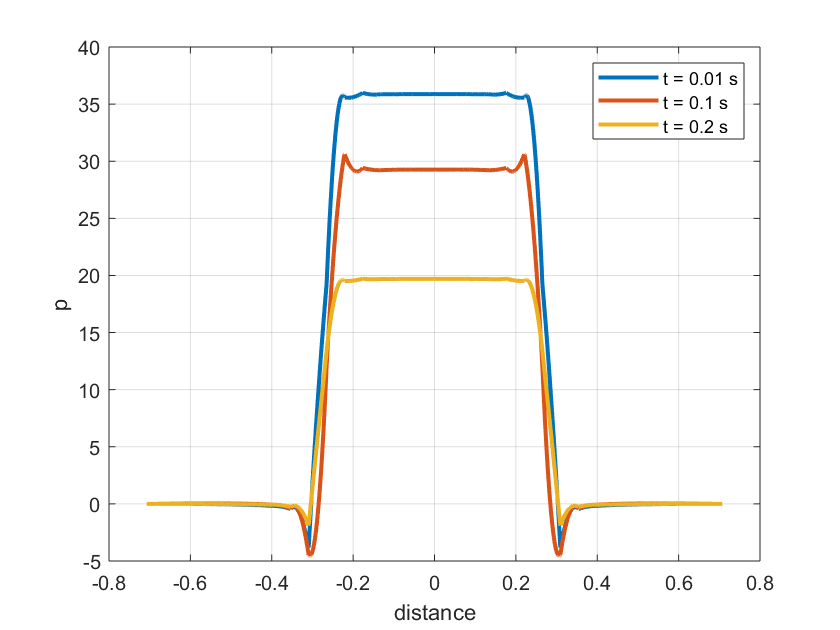}
        \caption{Level 2}
    \end{subfigure}%
    \begin{subfigure}[b]{0.32\textwidth}
        \includegraphics[width=\textwidth]{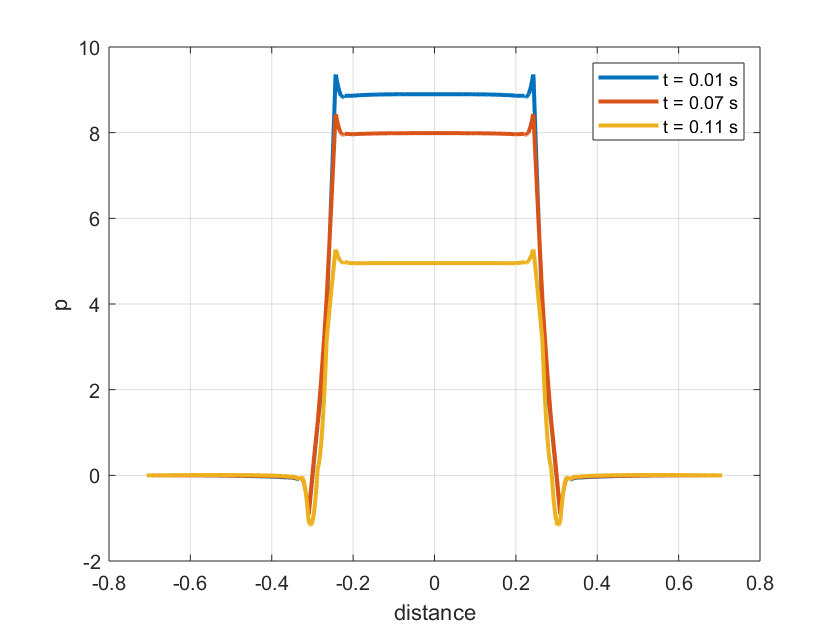}
        \caption{Level 3}
    \end{subfigure}
    \caption{Pressure $p$ over the whole domain at successive times for the three
    refinement levels of the elliptically displaced ring.}
    \label{fig:ellipseplot}
\end{figure}

As shown in Section~\ref{sec:intro} (Fig.~\ref{fig:fail}), the unstabilized
DLM-FE-IBM~\eqref{eqn:fdibm} produces a spurious collapse of the ring for this
problem. With the mixed formulation, incompressibility of the solid is enforced
and the ring returns to its circular equilibrium. Table~\ref{table:ellipse}
lists three refinement levels, all with solid-to-fluid mesh ratio
$h_s/h_x\sim 2$, and Fig.~\ref{fig:ellipse} shows the corresponding evolution.
For every level the initially elliptical ring relaxes to the circular
equilibrium configuration, confirming that the volumetric stabilization removes
the failure of Fig.~\ref{fig:fail}. Figure~\ref{fig:ellipseplot} reports the
pressure over the domain at successive times: the pressure inside the ring
decreases as the ring returns to equilibrium, with some pressure artefacts at the
solid--fluid interface.

\begin{table}[t]
    \caption{Mesh and refinement levels for the elliptically displaced ring.}
    \label{table:ellipse}
    \centering
\begin{tabularx}{0.85\textwidth}{@{}lYYYYY@{}} \toprule
    {} & {$h_x$} & {Fluid cells} & {Solid cells} & {Total DoFs} & {$h_s/h_x$} \\ \midrule
    Level 1 & \sfrac{1}{16} & 64 & 5 & 76 & $\sim 2$ \\
    Level 2 & \sfrac{1}{32} & 256 & 20 & 2776 & $\sim 2$ \\
    Level 3 & \sfrac{1}{64} & 1024 & 80 & 10592 & $\sim 2$ \\ \bottomrule
\end{tabularx}
\end{table}

\subsection{Radially stretched thick ring}
\label{sec:stretch}

Next we place a radially stretched thick ring at the centre of the fluid domain
and let it evolve. In the absence of fluid the ring would return to its reference
configuration; because both the fluid and the solid are incompressible, however,
the ring cannot contract, and a pressure difference must develop between the
interior and the exterior of the ring. This test therefore probes whether the
mixed formulation enforces incompressibility of the solid strongly enough to
sustain a physical pressure jump. The domain and ring geometry are as in
Section~\ref{sec:ellipse}; the ring is stretched radially by~$10\%$ with
$\mu_s = 10^4\ \mathrm{N/mm^2}$, $\rho_f = 1\ \mathrm{kg/mm^3}$ and
$\mu_f = 1\ \mathrm{N\cdot s/mm^2}$, and the system evolves for~$1$~s.

\begin{figure}[t]
    \centering
    \begin{subfigure}[b]{0.32\textwidth}
        \includegraphics[width=\textwidth]{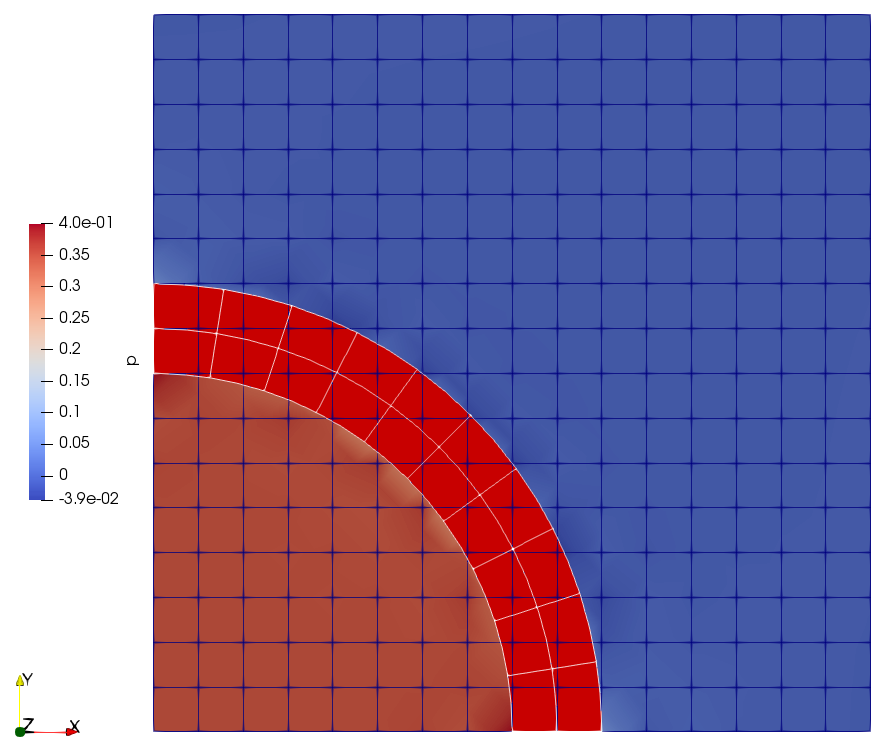}
        \caption{Case A ($h_s/h_x\sim 2$)}
    \end{subfigure}%
    \begin{subfigure}[b]{0.32\textwidth}
        \includegraphics[width=\textwidth]{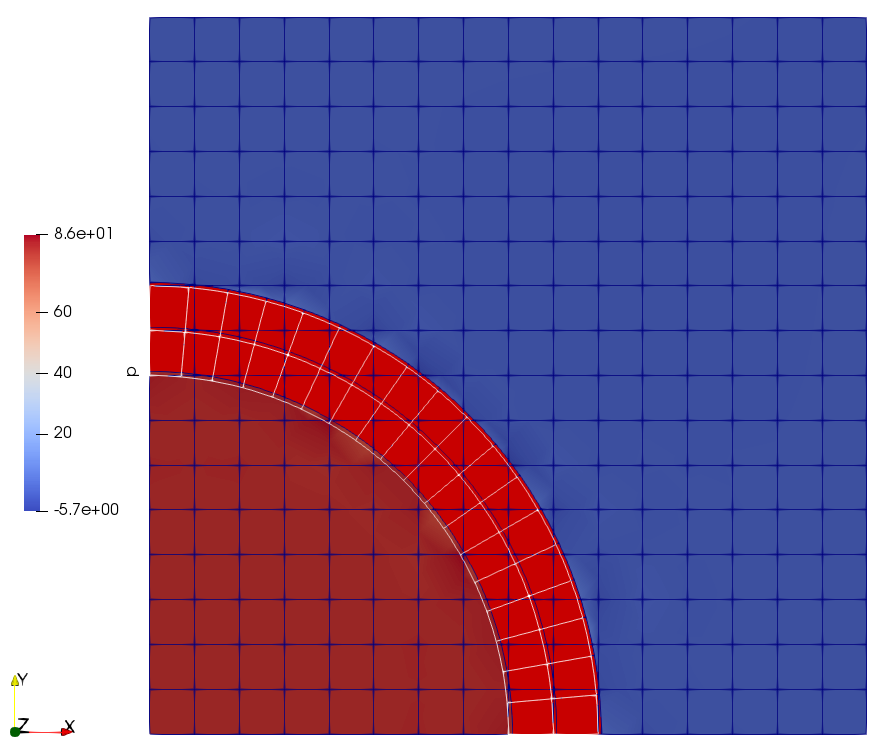}
        \caption{Case B ($h_s/h_x\sim 1$)}
    \end{subfigure}%
    \begin{subfigure}[b]{0.32\textwidth}
        \includegraphics[width=\textwidth]{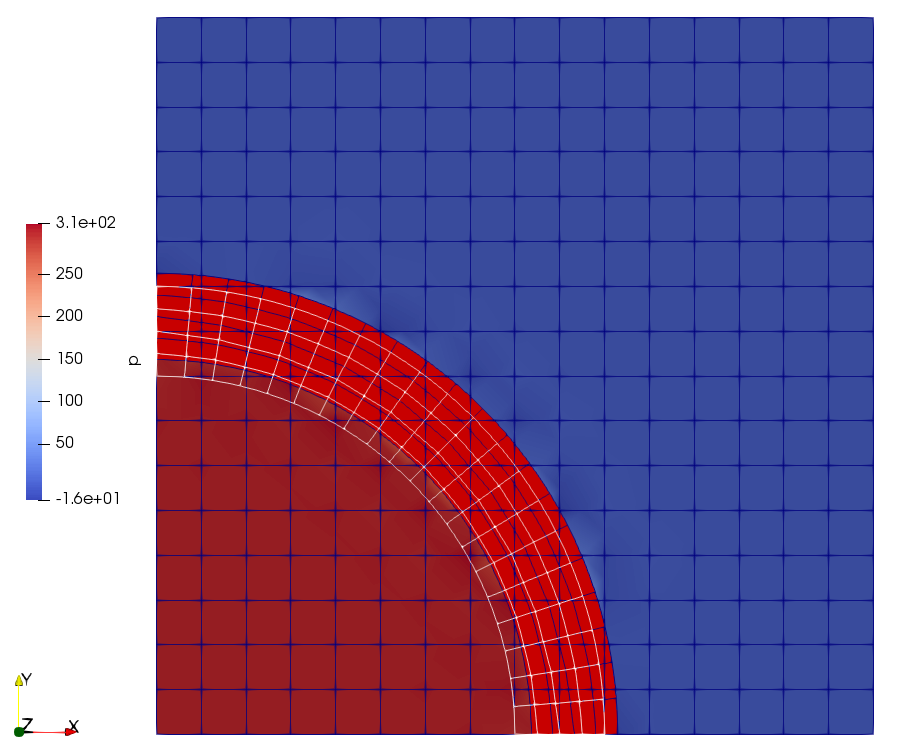}
        \caption{Case C ($h_s/h_x\sim\sfrac{1}{2}$)}
    \end{subfigure}
    \caption{Radially stretched ring at $t = 0.99$~s for different solid-to-fluid
    mesh ratios (cases A--C of Table~\ref{table:stretch}). Only for
    $h_s/h_x\sim\tfrac{1}{2}$ (case~C) does the incompressible ring retain its
    stretched configuration; for larger ratios the ring incorrectly relaxes to
    its reference shape (shown as the white wireframe).}
    \label{fig:stretch1}
\end{figure}

\begin{figure}[t]
    \centering
    \begin{subfigure}[b]{0.32\textwidth}
        \includegraphics[width=\textwidth]{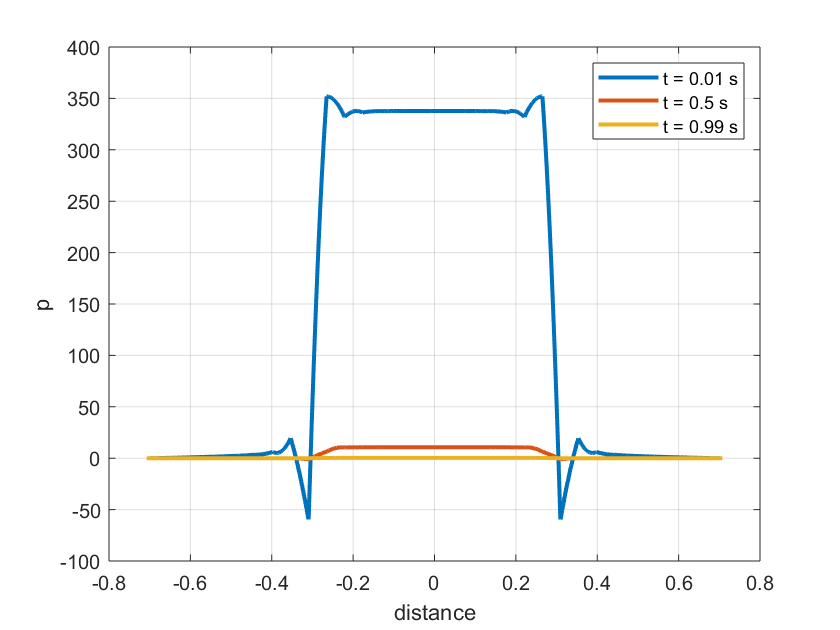}
        \caption{Case A}
    \end{subfigure}%
    \begin{subfigure}[b]{0.32\textwidth}
        \includegraphics[width=\textwidth]{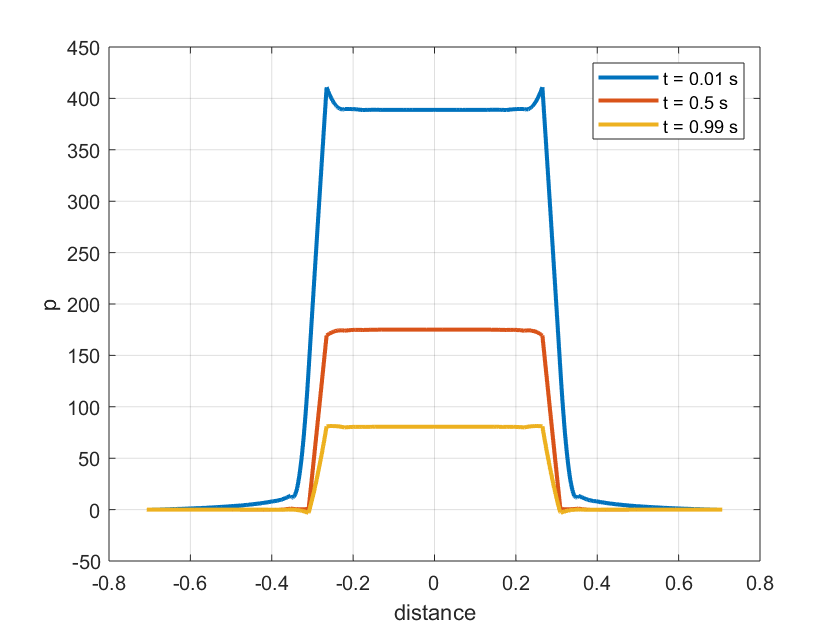}
        \caption{Case B}
    \end{subfigure}%
    \begin{subfigure}[b]{0.32\textwidth}
        \includegraphics[width=\textwidth]{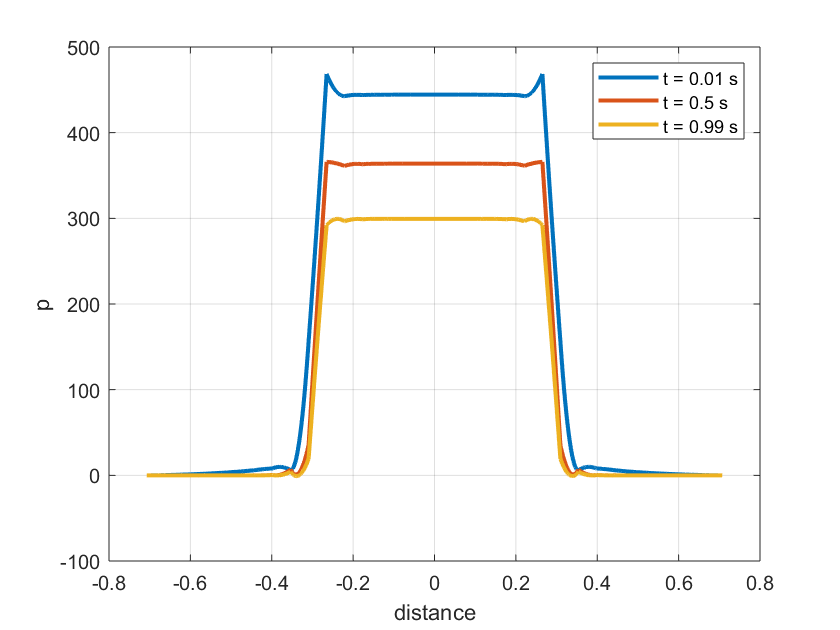}
        \caption{Case C}
    \end{subfigure}
    \caption{Pressure difference between the interior and exterior of the ring
    over time for cases A--C. The difference decays for cases~A and~B but is
    sustained for case~C.}
    \label{fig:stretchplot1}
\end{figure}

The cases in Table~\ref{table:stretch} vary both the solid-to-fluid mesh ratio
$h_s/h_x$ (cases A, B, C) and the overall refinement (cases C, D, E).
Figures~\ref{fig:stretch1} and~\ref{fig:stretchplot1} show that the ratio
$h_s/h_x$ is decisive. For $h_s/h_x\sim 2$ (case~A) and $\sim 1$ (case~B) the
solid relaxes back toward its reference configuration and the pressure difference
decays to zero---i.e.\ incompressibility of the solid is not maintained---whereas
for $h_s/h_x\sim\tfrac{1}{2}$ (case~C) the ring holds its stretched
configuration and a significant, sustained pressure difference develops, as
physically expected.

\begin{table}[t]
    \caption{Mesh and refinement cases for the radially stretched ring.}
    \label{table:stretch}
    \centering
\begin{tabularx}{0.85\textwidth}{@{}lYYYYY@{}} \toprule
    {} & {$h_x$} & {Fluid cells} & {Solid cells} & {Total DoFs} & {$h_s/h_x$} \\ \midrule
    Case A & \sfrac{1}{32} & 256 & 20 & 2776 & $\sim 2$ \\
    Case B & \sfrac{1}{32} & 256 & 36 & 3008 & $\sim 1$ \\
    Case C & \sfrac{1}{32} & 256 & 60 & 3520 & $\sim\sfrac{1}{2}$ \\
    Case D & \sfrac{1}{16} & 64 & 20 & 968 & $\sim\sfrac{1}{2}$ \\
    Case E & \sfrac{1}{64} & 1024 & 320 & 13400 & $\sim\sfrac{1}{2}$ \\ \bottomrule
\end{tabularx}
\end{table}

\begin{figure}[t]
    \centering
    \begin{subfigure}[b]{0.32\textwidth}
        \includegraphics[width=\textwidth]{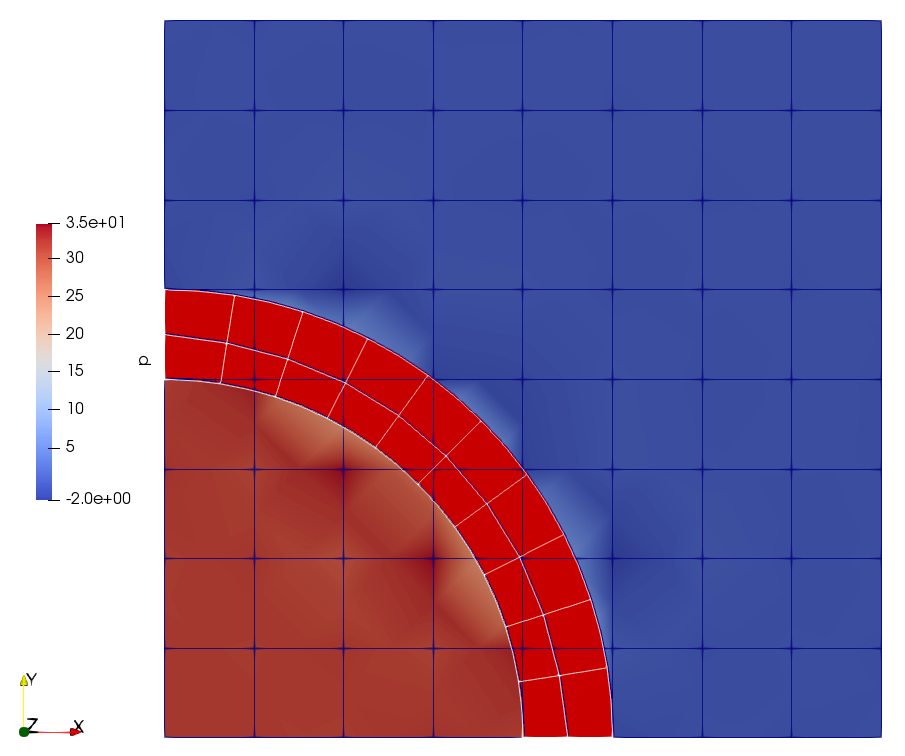}
        \caption{Case D (coarse)}
    \end{subfigure}%
    \begin{subfigure}[b]{0.32\textwidth}
        \includegraphics[width=\textwidth]{stretch2to1_fin_1.png}
        \caption{Case C (medium)}
    \end{subfigure}%
    \begin{subfigure}[b]{0.32\textwidth}
        \includegraphics[width=\textwidth]{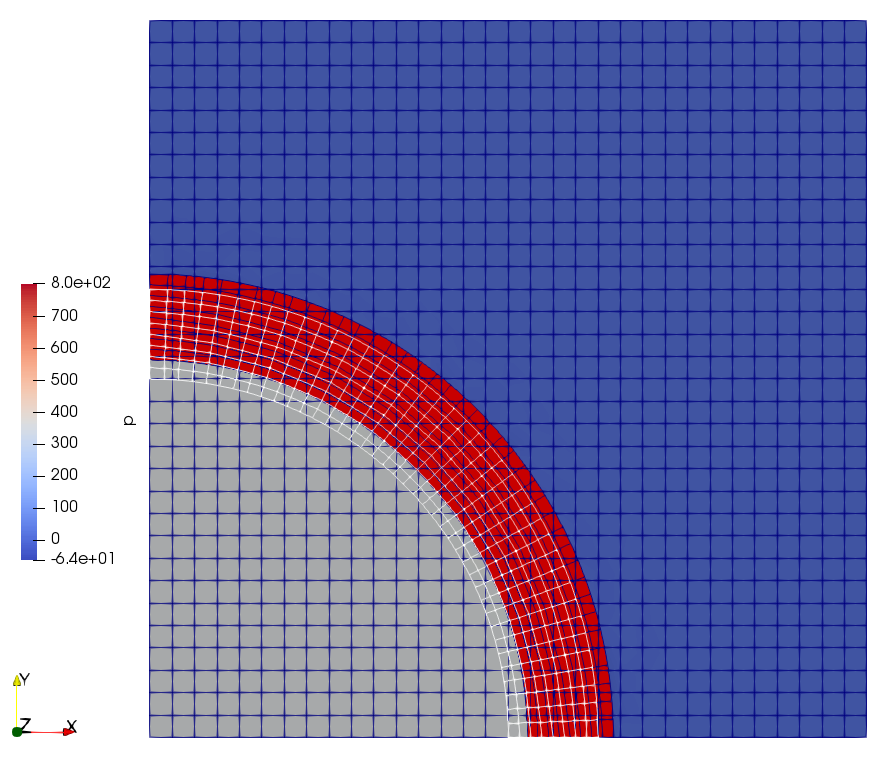}
        \caption{Case E (fine)}
    \end{subfigure}
    \caption{Radially stretched ring at $t = 0.99$~s for increasing refinement at
    fixed $h_s/h_x\sim\tfrac{1}{2}$ (cases D, C, E of
    Table~\ref{table:stretch}). The stretched configuration is recovered as the
    mesh is refined.}
    \label{fig:stretch2}
\end{figure}

\begin{figure}[t]
    \centering
    \begin{subfigure}[b]{0.32\textwidth}
        \includegraphics[width=\textwidth]{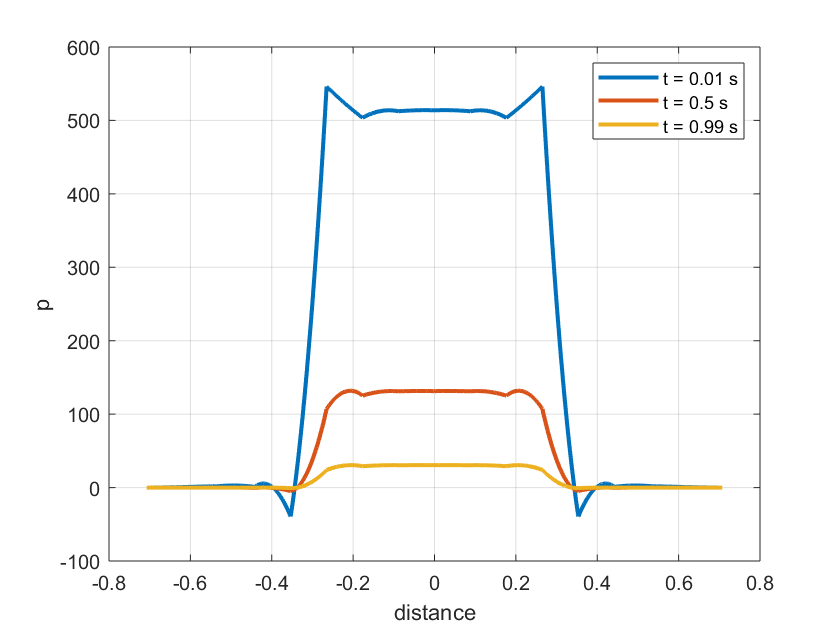}
        \caption{Case D}
    \end{subfigure}%
    \begin{subfigure}[b]{0.32\textwidth}
        \includegraphics[width=\textwidth]{stretch2to1_1.png}
        \caption{Case C}
    \end{subfigure}%
    \begin{subfigure}[b]{0.32\textwidth}
        \includegraphics[width=\textwidth]{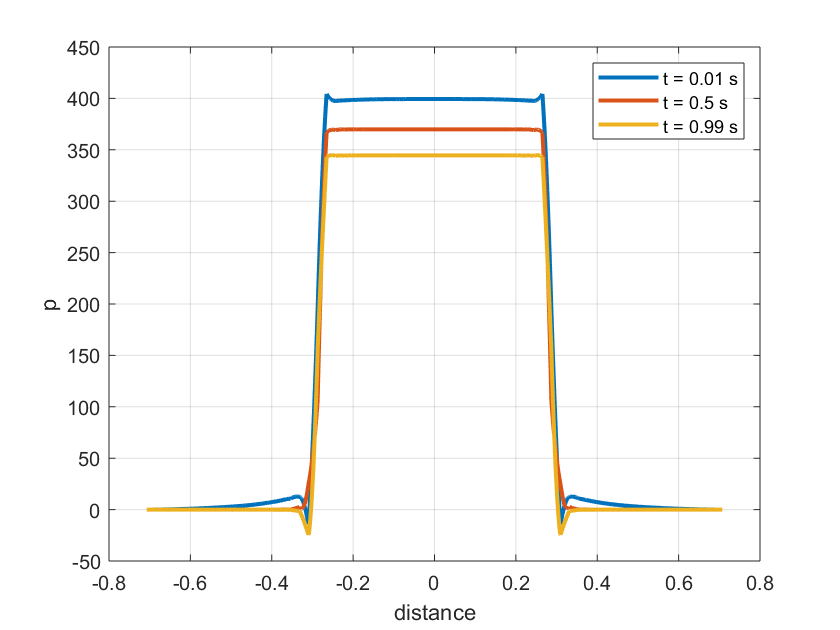}
        \caption{Case E}
    \end{subfigure}
    \caption{Pressure difference over time for the refinement study (cases D, C,
    E). Convergence toward a sustained pressure difference is observed as the mesh
    is refined.}
    \label{fig:stretchplot2}
\end{figure}

Fixing $h_s/h_x\sim\tfrac{1}{2}$ and varying the overall refinement (cases D, C,
E) shows the expected mesh convergence: Figs.~\ref{fig:stretch2}
and~\ref{fig:stretchplot2} indicate that the coarsest mesh (case~D) again fails
to sustain the pressure difference, whereas refining the mesh (case~E) yields the
expected behaviour, with the final-time pressure difference remaining close to
its initial value and the interface pressure artefacts diminishing under
refinement. A small amount of residual pressure leakage remains and is expected
to vanish with further refinement.

\begin{figure}[t]
    \centering
    \includegraphics[width=0.6\textwidth]{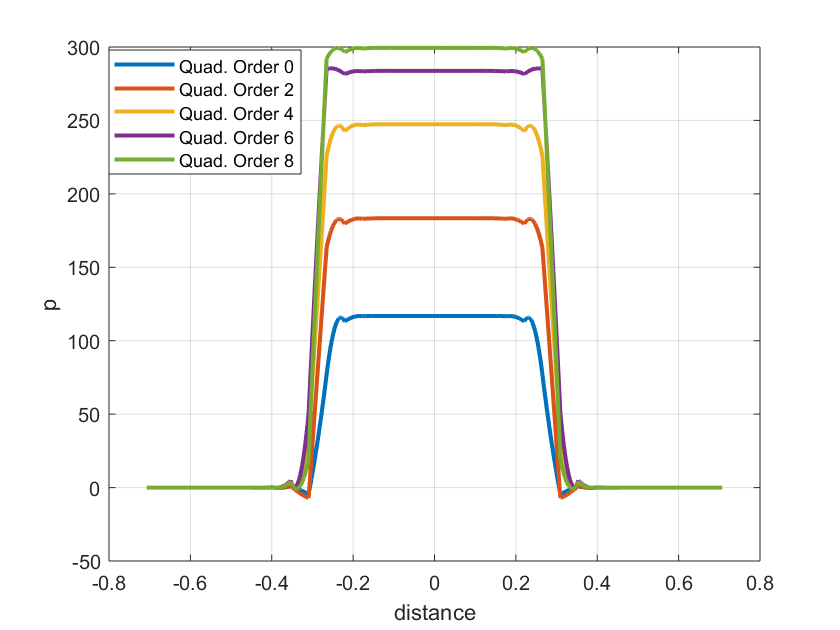}
    \caption{Effect of the extra quadrature order on the domain pressure for
    case~C of the radially stretched ring.}
    \label{fig:stretch_quad}
\end{figure}

Finally, Fig.~\ref{fig:stretch_quad} examines the effect of the extra quadrature
order supplied during assembly of the non-matching coupling terms, for case~C.
As the extra quadrature order increases from~$0$ to~$8$, the sustained pressure
difference increases and approaches the value at the initial time, confirming
that a sufficient quadrature order is needed to integrate the non-matching
coupling accurately; order~$8$ was used for the results above.

\subsection{Disk falling under gravity in a viscous fluid}
\label{sec:falling}

\begin{figure}[t]
    \centering
    \includegraphics[width=0.34\linewidth]{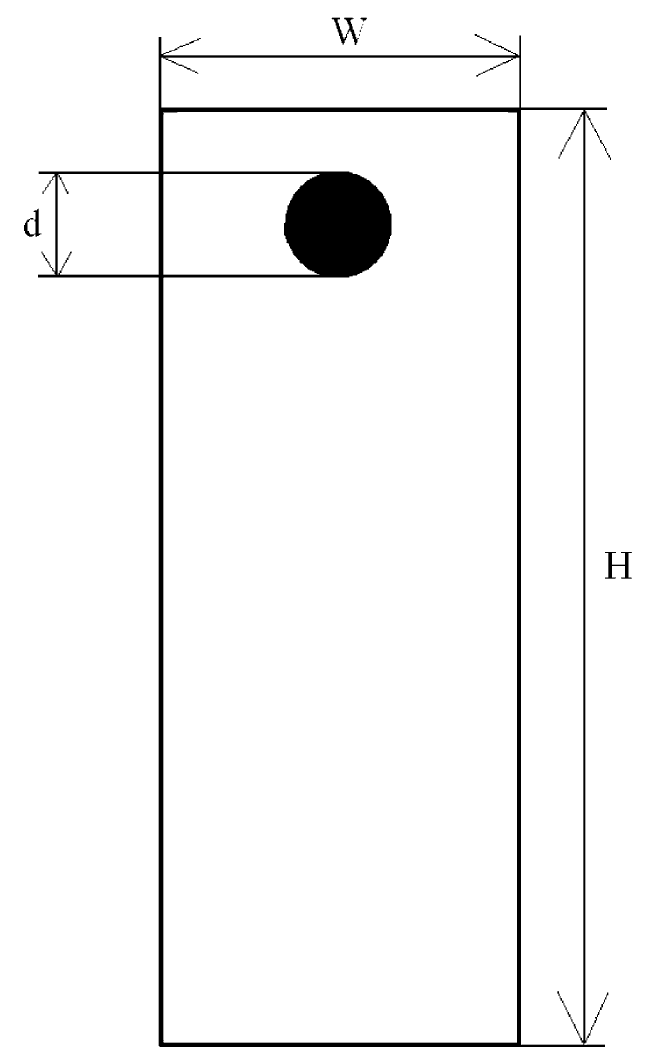}
    \caption{Geometry of the disk falling in a 2D control volume.}
    \label{fig:drop}
\end{figure}

To test the method against a quantitative reference we study a disk (the
mid-plane of a cylinder) falling under gravity in an incompressible viscous
fluid, following~\cite{roy,zhang}. A disk of diameter $d$ is released from rest
in a rectangular control volume of height $H$ and width $W$
(Fig.~\ref{fig:drop}); as it sinks we track the centre of the disk and its
velocity~$v$, which approaches a terminal velocity~$v_T$. The fixed parameters
are $H = 2.0$~cm, $W = 1.0$~cm, $d = 0.1$~cm, $\rho_f = 1.0\ \mathrm{g/cm^3}$,
$\mu_s = 10^3\ \mathrm{dyn/cm^2}$ and $g = 981\ \mathrm{cm/s^2}$; the control
volume uses $16{,}000$ cells, the solid $40$ cells, for $116{,}844$ total
degrees of freedom.

We compare the computed terminal velocity with the empirical solution for a
rigid disk settling in a viscous fluid~\cite{clift},
\begin{equation}\label{eqn:exact}
v_E = \frac{(\rho_s-\rho_f)\,g\,r^2}{4\mu_f}
   \left(\ln\!\Bigl(\tfrac{L}{r}\Bigr) - 0.9157
   + 1.7244\Bigl(\tfrac{r}{L}\Bigr)^2
   - 1.7302\Bigl(\tfrac{r}{L}\Bigr)^4\right),
\end{equation}
where $r = d/2$ and $L = W/2$. Table~\ref{table:drop1} varies the solid density
at fixed viscosity $\mu_f = 1.0$~P, and Table~\ref{table:drop2} varies the fluid
viscosity at fixed density $\rho_s = 4.0\ \mathrm{g/cm^3}$. In all cases the
computed terminal velocity $v_T$ agrees with the empirical value $v_E$ to within
about~$1\%$.

\begin{table}[t]
    \caption{Terminal velocity for varying solid density (constant
    $\mu_f = 1.0$~P): computed $v_T$ versus the empirical solution
    $v_E$ of~\eqref{eqn:exact}.}
    \label{table:drop1}
    \centering
\begin{tabularx}{0.7\textwidth}{@{}YYYY@{}} \toprule
    {$\rho_s\ (\mathrm{g/cm^3})$} & {$\mu_f\ (\mathrm{P})$} & {$v_T\ (\mathrm{cm/s})$} & {$v_E\ (\mathrm{cm/s})$} \\ \midrule
    2.0 & 1.0 & 0.866 & 0.8608 \\
    3.0 & 1.0 & 1.722 & 1.7216 \\
    4.0 & 1.0 & 2.566 & 2.5824 \\ \bottomrule
\end{tabularx}
\end{table}

\begin{table}[t]
    \caption{Terminal velocity for varying fluid viscosity (constant
    $\rho_s = 4.0\ \mathrm{g/cm^3}$): computed $v_T$ versus the empirical
    solution $v_E$.}
    \label{table:drop2}
    \centering
\begin{tabularx}{0.7\textwidth}{@{}YYYY@{}} \toprule
    {$\rho_s\ (\mathrm{g/cm^3})$} & {$\mu_f\ (\mathrm{P})$} & {$v_T\ (\mathrm{cm/s})$} & {$v_E\ (\mathrm{cm/s})$} \\ \midrule
    4.0 & 1.0 & 2.566 & 2.5824 \\
    4.0 & 2.0 & 1.297 & 1.2912 \\
    4.0 & 3.0 & 0.866 & 0.8608 \\ \bottomrule
\end{tabularx}
\end{table}

\begin{figure}[t]
    \centering
    \begin{subfigure}[b]{0.24\textwidth}
        \includegraphics[width=\textwidth]{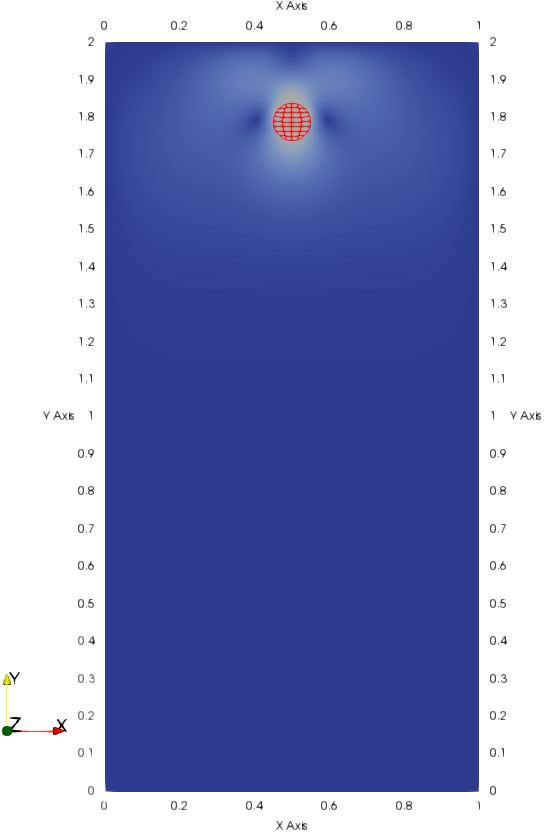}
        \caption{$t = 0$}
    \end{subfigure}%
    \begin{subfigure}[b]{0.24\textwidth}
        \includegraphics[width=\textwidth]{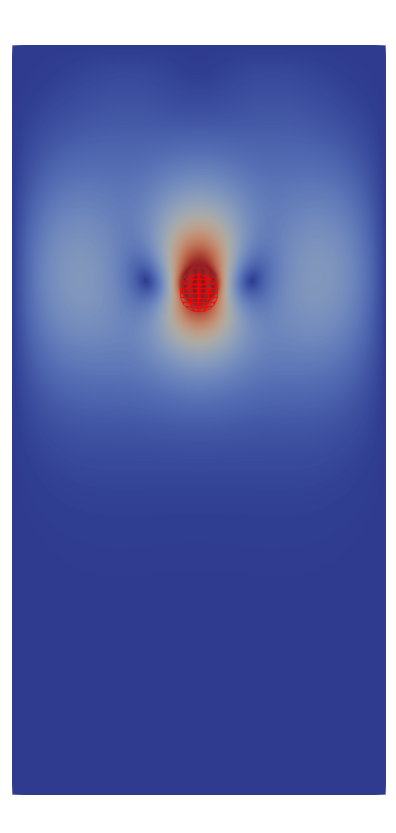}
        \caption{$t = 0.2$}
    \end{subfigure}%
    \begin{subfigure}[b]{0.24\textwidth}
        \includegraphics[width=\textwidth]{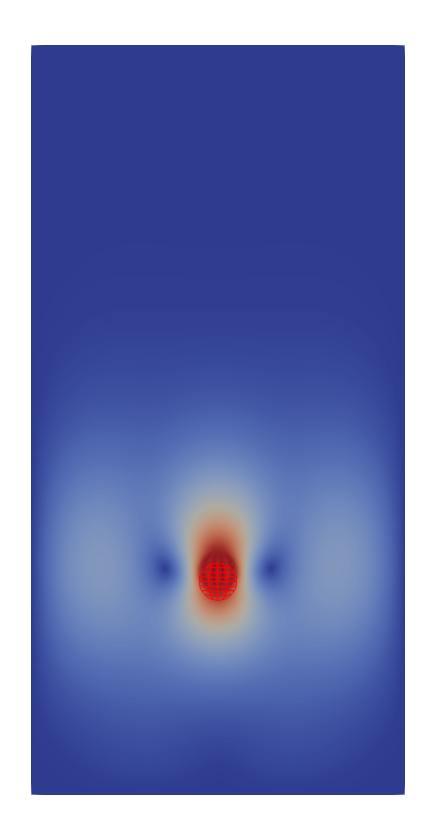}
        \caption{$t = 0.5$}
    \end{subfigure}%
    \begin{subfigure}[b]{0.24\textwidth}
        \includegraphics[width=\textwidth]{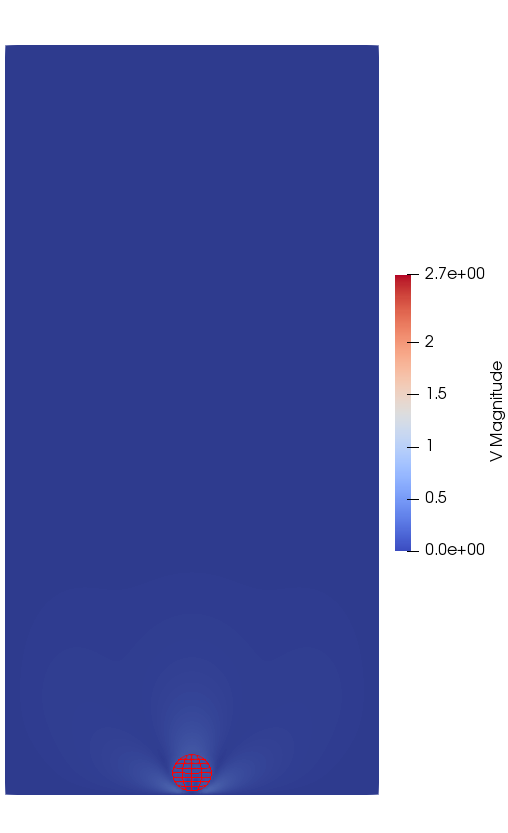}
        \caption{$t = 0.8$}
    \end{subfigure}
    \caption{Fluid velocity field and disk position at successive times for
    $\rho_s = 4.0\ \mathrm{g/cm^3}$, $\mu_f = 1.0$~P.}
    \label{fig:fallingball}
\end{figure}

\begin{figure}[t]
    \centering
    \begin{subfigure}[b]{0.49\textwidth}
        \includegraphics[width=\textwidth]{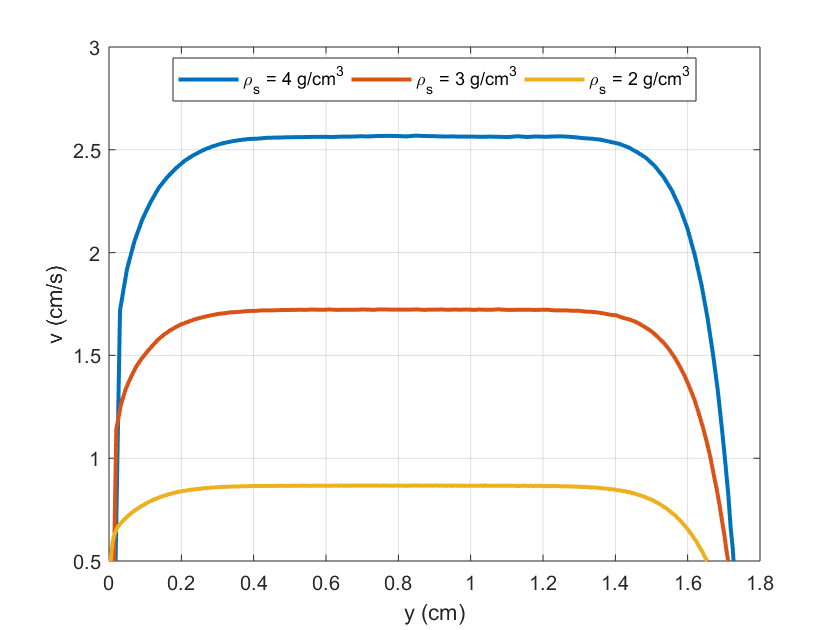}
        \caption{Disk velocity vs.\ axial position}
    \end{subfigure}%
    \begin{subfigure}[b]{0.49\textwidth}
        \includegraphics[width=\textwidth]{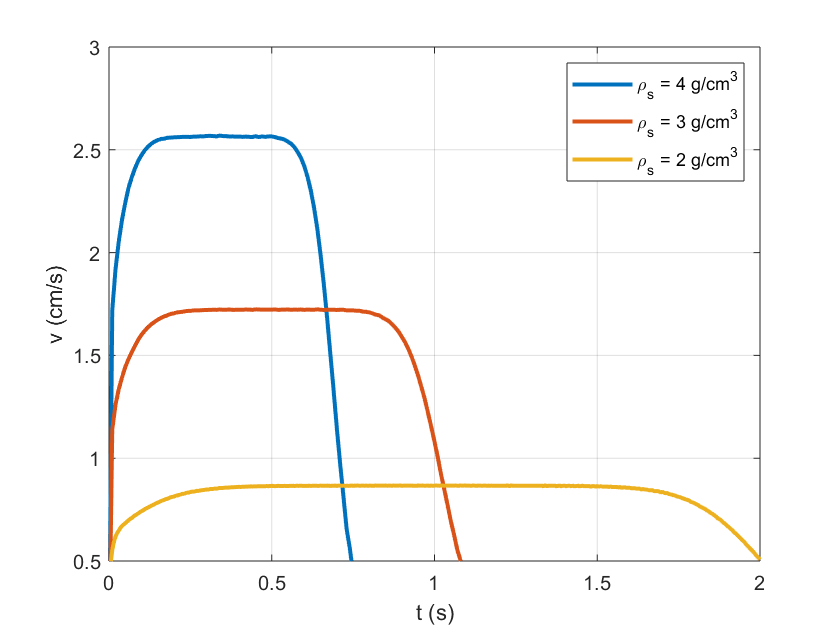}
        \caption{Disk velocity vs.\ time}
    \end{subfigure}
    \caption{Effect of solid density on the disk velocity.}
    \label{fig:drop1}
\end{figure}

\begin{figure}[t]
    \centering
    \begin{subfigure}[b]{0.49\textwidth}
        \includegraphics[width=\textwidth]{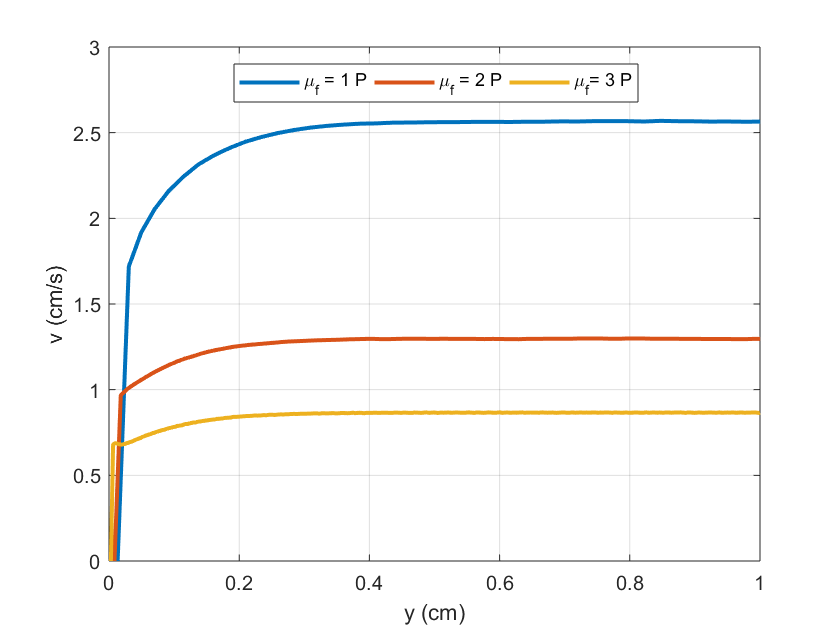}
        \caption{Disk velocity vs.\ axial position}
    \end{subfigure}%
    \begin{subfigure}[b]{0.49\textwidth}
        \includegraphics[width=\textwidth]{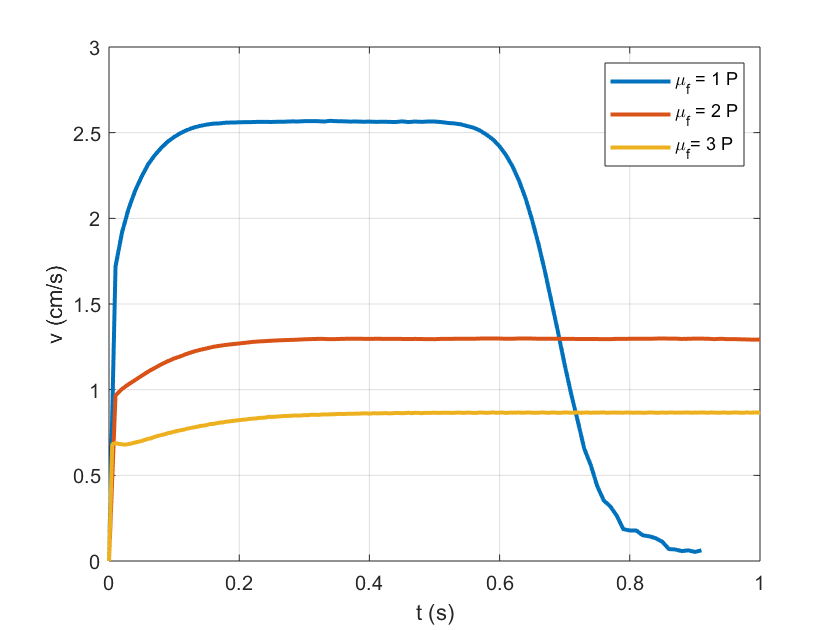}
        \caption{Disk velocity vs.\ time}
    \end{subfigure}
    \caption{Effect of fluid viscosity on the disk velocity.}
    \label{fig:drop2}
\end{figure}

\begin{figure}[t]
    \centering
    \begin{subfigure}[b]{0.49\textwidth}
        \includegraphics[width=\textwidth]{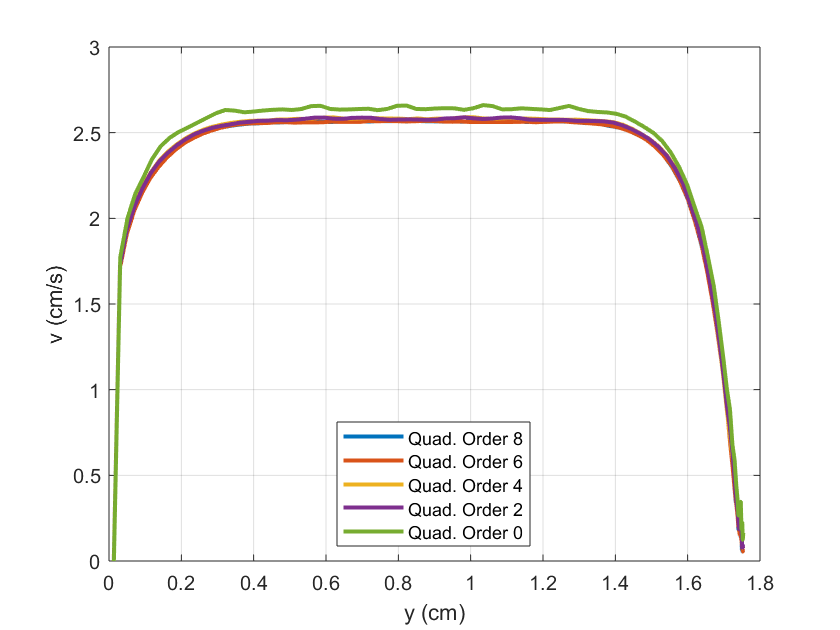}
        \caption{Disk velocity vs.\ axial position}
    \end{subfigure}%
    \begin{subfigure}[b]{0.49\textwidth}
        \includegraphics[width=\textwidth]{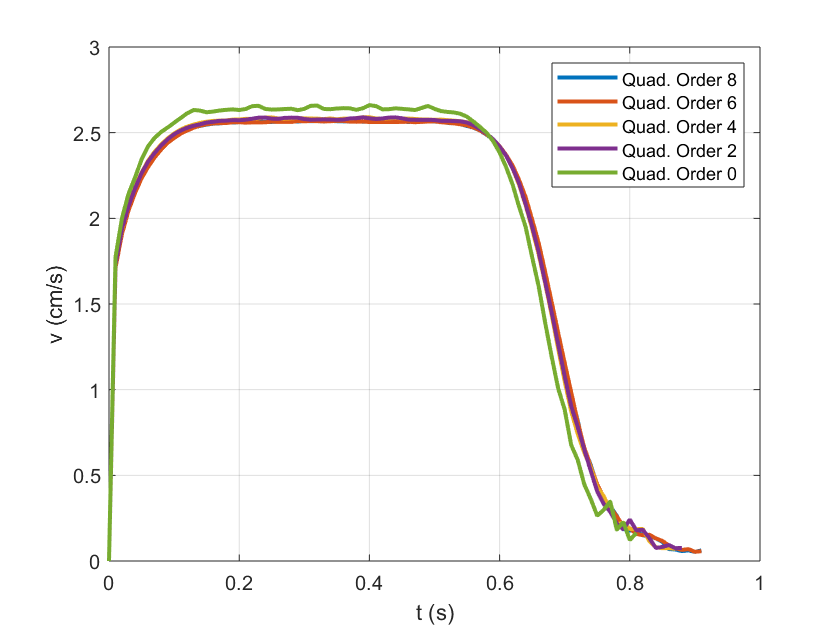}
        \caption{Disk velocity vs.\ time}
    \end{subfigure}
    \caption{Effect of the extra quadrature order on the disk velocity for
    $\rho_s = 4.0\ \mathrm{g/cm^3}$, $\mu_f = 1.0$~P. Orders~4, 6 and~8 coincide.}
    \label{fig:dropquad}
\end{figure}

Figure~\ref{fig:fallingball} shows the velocity field and disk position at
successive times, and Figs.~\ref{fig:drop1} and~\ref{fig:drop2} show that the
disk accelerates from rest, reaches a terminal velocity over the central portion
of the domain, and decelerates as it nears the bottom wall. The terminal velocity
increases with solid density and decreases with fluid viscosity, as expected and
in agreement with~\cite{roy}. Figure~\ref{fig:dropquad} reports the sensitivity
to the extra quadrature order for $\rho_s = 4.0\ \mathrm{g/cm^3}$,
$\mu_f = 1.0$~P: the terminal velocity converges as the order increases, and
orders~$4$, $6$ and~$8$ are essentially indistinguishable, so that extra
quadrature of order~$4$ already suffices for this problem. We emphasize that,
unlike an ALE simulation of the same problem, which would require periodic
remeshing because of the conforming fluid--solid interface, the present immersed
computation keeps the mesh topology fixed and requires no remeshing.

\subsection{An open challenge: the inflating ring}
\label{sec:inflate}

The three benchmarks above are reproduced correctly by the mixed formulation.
For completeness, and in the interest of transparency, we also report a more
demanding pressure-loaded case in which the method does not yet perform robustly.
We take the same ring as in Sections~\ref{sec:ellipse} and~\ref{sec:stretch},
placed at the centre of the domain, and inject a volume of fluid at the centre of
the ring over an initial interval. As the fluid is injected the ring should
stretch and thin; once injection stops, the incompressibility of the fluid and
the solid should prevent the ring from relaxing, sustaining a pressure difference
across the ring wall. The material and domain parameters are as before, and the
refinement cases are listed in Table~\ref{table:inflate}.

\begin{table}[t]
    \caption{Mesh and refinement cases for the inflating ring.}
    \label{table:inflate}
    \centering
\begin{tabularx}{0.85\textwidth}{@{}lYYYYY@{}} \toprule
    {} & {$h_x$} & {Fluid cells} & {Solid cells} & {Total DoFs} & {$h_s/h_x$} \\ \midrule
    Case I   & \sfrac{1}{32} & 256 & 20 & 2776 & $\sim 2$ \\
    Case II  & \sfrac{1}{32} & 256 & 36 & 3008 & $\sim 1$ \\
    Case III & \sfrac{1}{64} & 1024 & 144 & 11408 & $\sim 1$ \\ \bottomrule
\end{tabularx}
\end{table}

\begin{figure}[t]
    \centering
    \begin{subfigure}[b]{0.32\textwidth}
        \includegraphics[width=\textwidth]{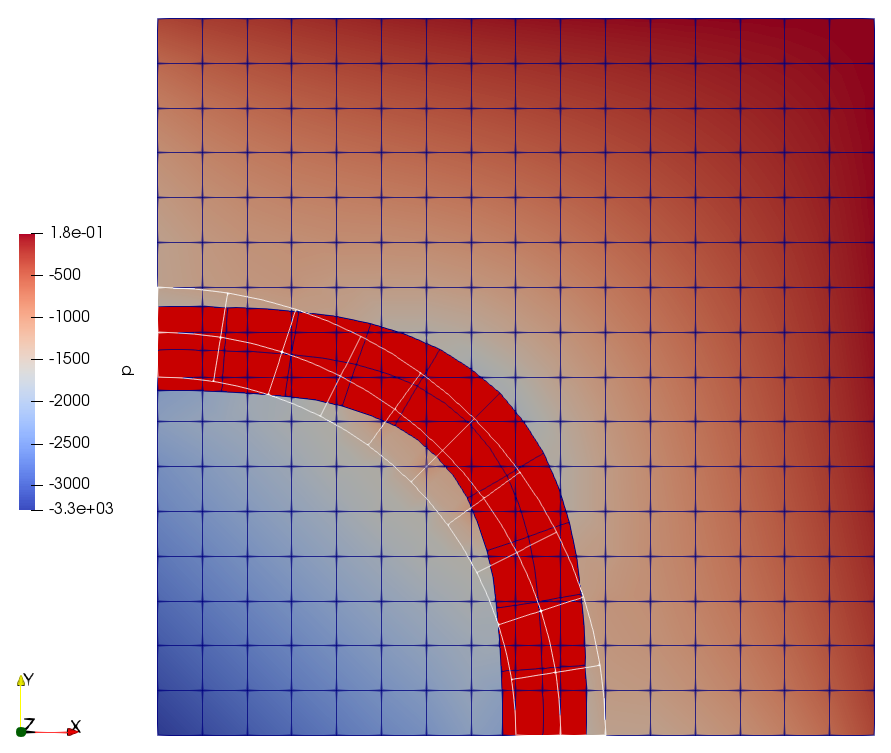}
        \caption{$t = 0.01$ s}
    \end{subfigure}%
    \begin{subfigure}[b]{0.32\textwidth}
        \includegraphics[width=\textwidth]{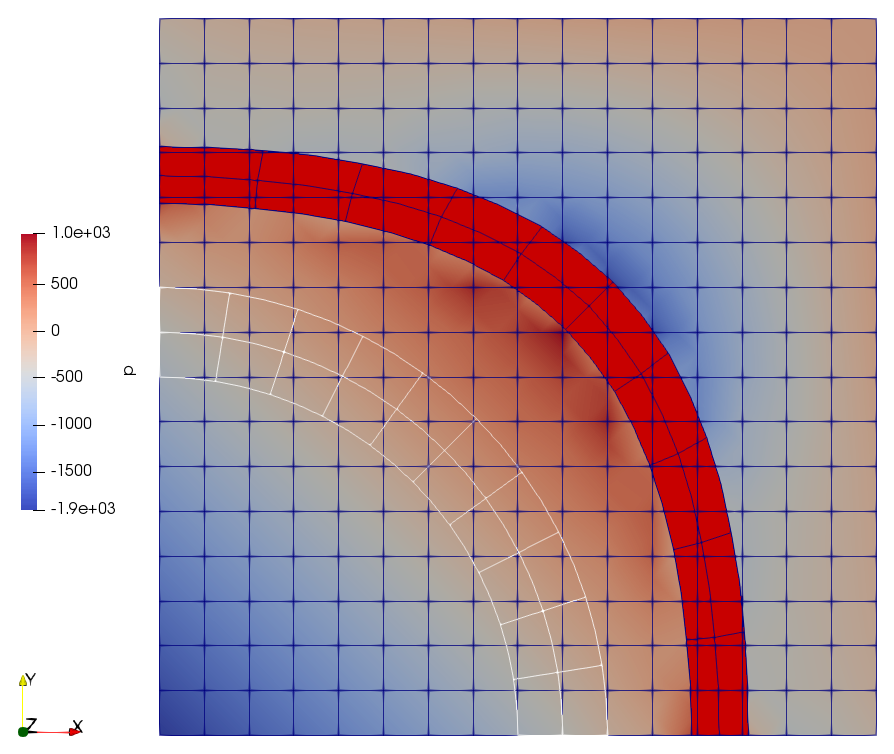}
        \caption{$t = 0.18$ s}
    \end{subfigure}%
    \begin{subfigure}[b]{0.32\textwidth}
        \includegraphics[width=\textwidth]{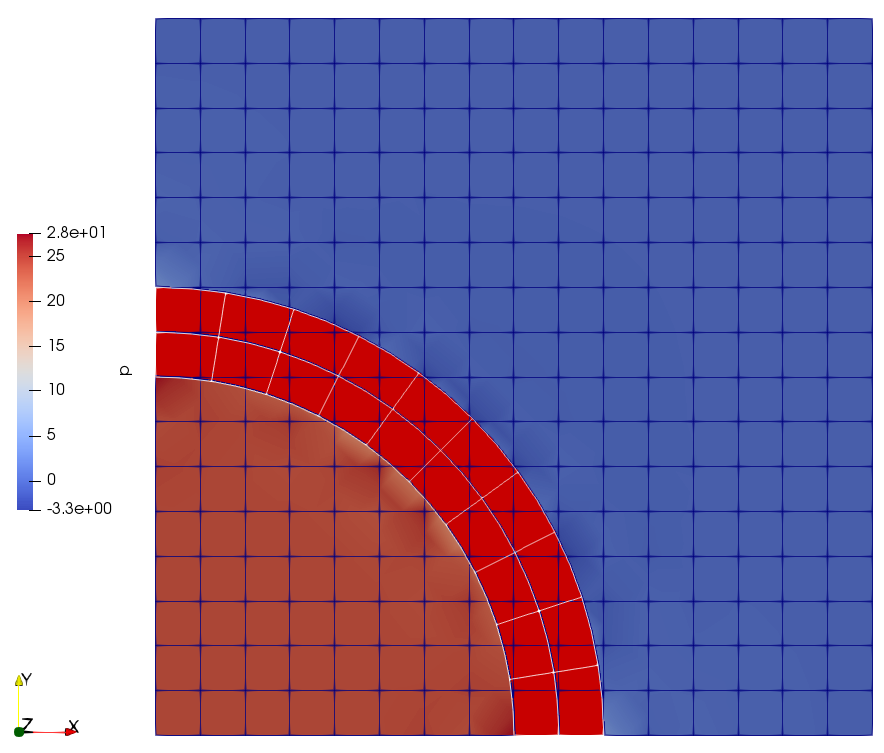}
        \caption{$t = 0.85$ s}
    \end{subfigure}
    \caption{Inflating ring, case~I ($h_s/h_x\sim 2$). The ring stretches during
    injection but then incorrectly relaxes to its reference configuration.}
    \label{fig:inflate1to2}
\end{figure}

\begin{figure}[t]
    \centering
    \begin{subfigure}[b]{0.32\textwidth}
        \includegraphics[width=\textwidth]{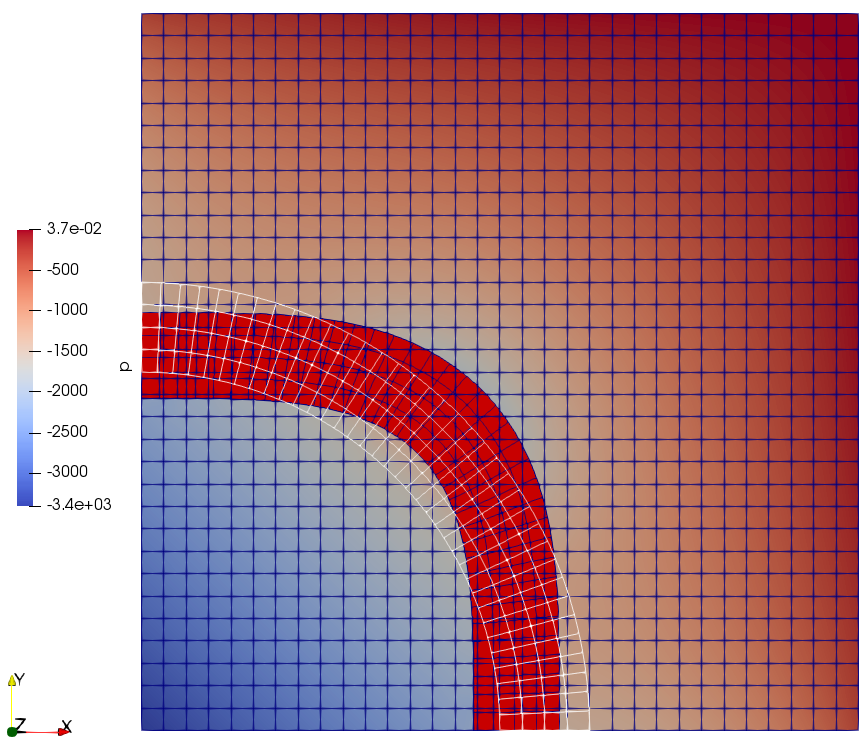}
        \caption{Case III: $t = 0.01$ s}
    \end{subfigure}%
    \begin{subfigure}[b]{0.32\textwidth}
        \includegraphics[width=\textwidth]{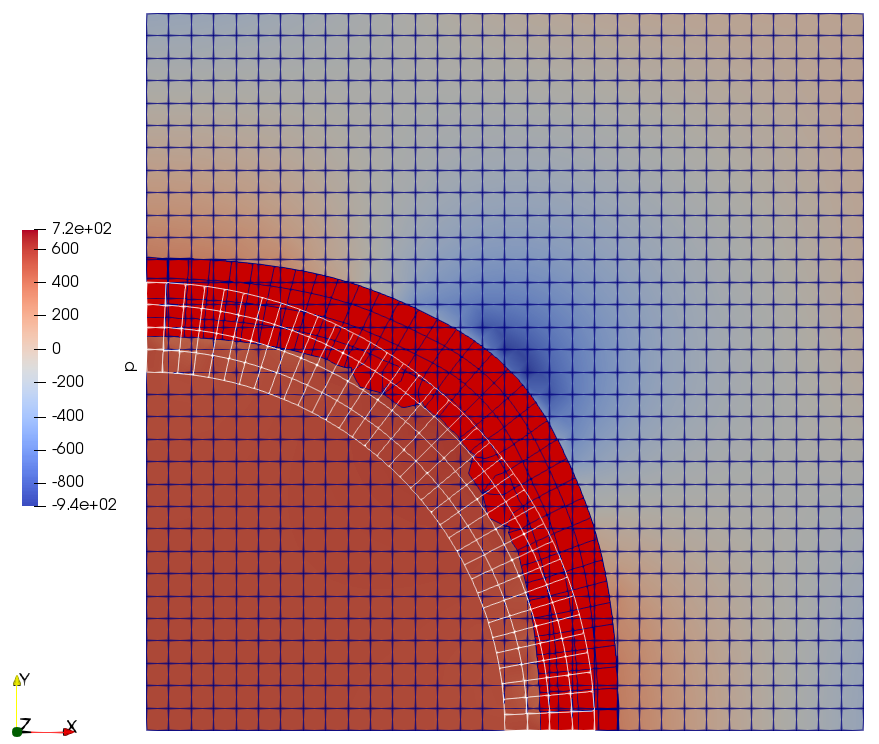}
        \caption{Case III: $t = 0.25$ s}
    \end{subfigure}%
    \begin{subfigure}[b]{0.32\textwidth}
        \includegraphics[width=\textwidth]{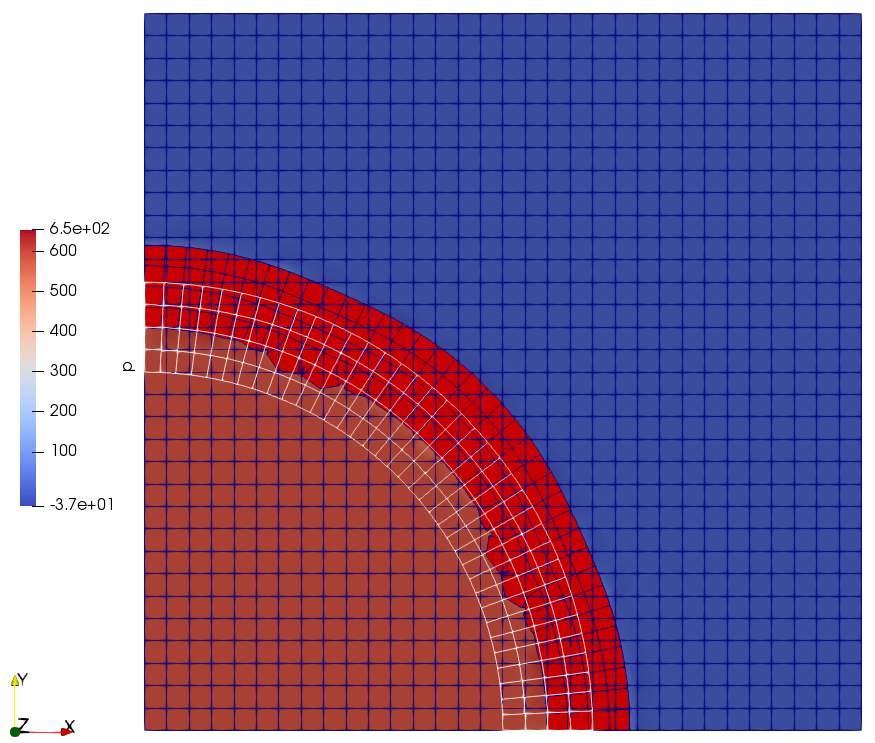}
        \caption{Case III: $t = 0.5$ s}
    \end{subfigure}
    \caption{Inflating ring, case~III ($h_s/h_x\sim 1$, refined). A pressure
    field develops inside the ring, but residual volumetric instabilities remain.}
    \label{fig:inflate1to1}
\end{figure}

\begin{figure}[t]
    \centering
    \begin{subfigure}[b]{0.49\textwidth}
        \includegraphics[width=\textwidth]{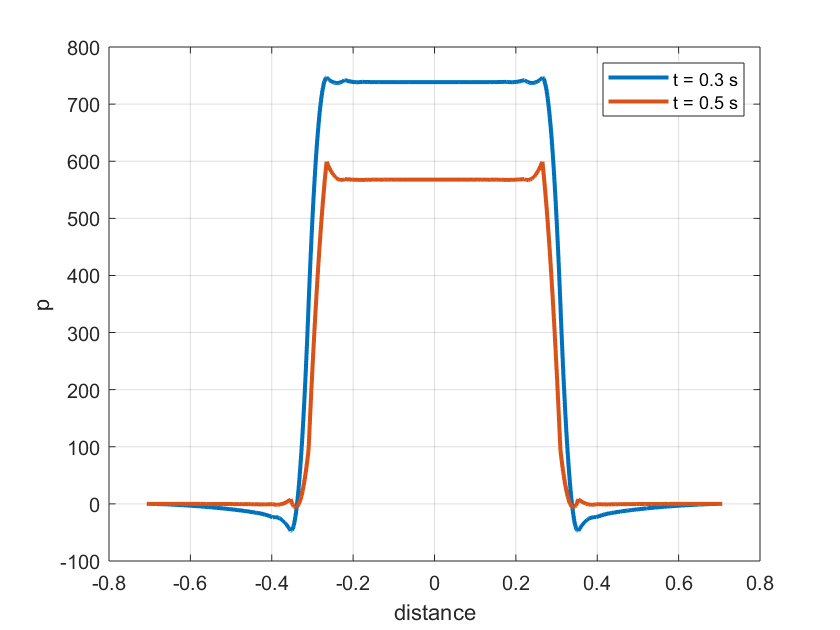}
        \caption{Case II}
    \end{subfigure}%
    \begin{subfigure}[b]{0.49\textwidth}
        \includegraphics[width=\textwidth]{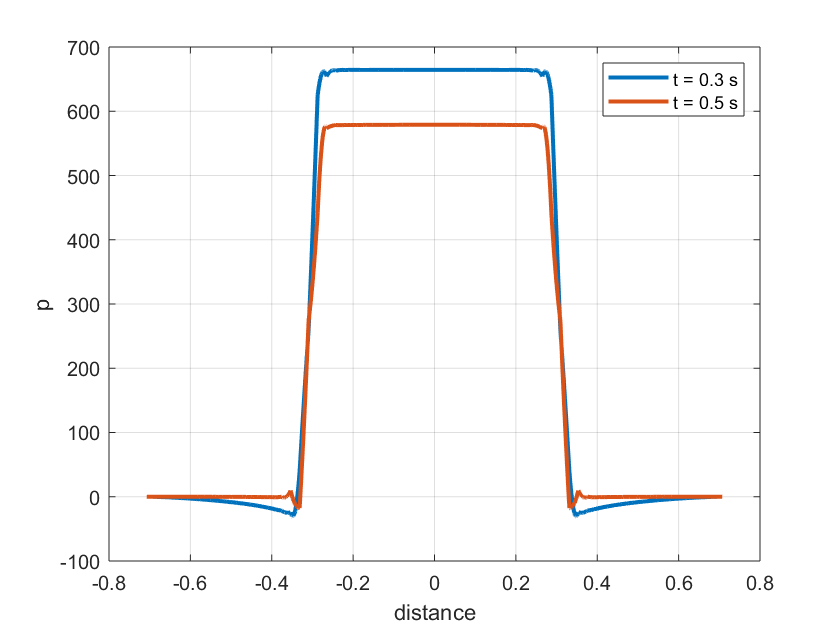}
        \caption{Case III}
    \end{subfigure}
    \caption{Pressure over the domain at successive times for the inflating ring,
    cases II and III. Refinement (case~III) yields a lower and smoother pressure
    difference.}
    \label{fig:inflateplot}
\end{figure}

For $h_s/h_x\sim 2$ (case~I) the ring stretches and thins as expected during
injection, but subsequently relaxes back to its reference configuration
(Fig.~\ref{fig:inflate1to2}), i.e.\ incompressibility of the solid is not
maintained after loading---the same trend seen for the stretched ring at large
$h_s/h_x$. Reducing the ratio to $h_s/h_x\sim 1$ (cases II and III) does develop
a pressure field inside the ring (Figs.~\ref{fig:inflate1to1}
and~\ref{fig:inflateplot}), and refinement (case~III) produces a lower and
smoother pressure difference than case~II; however, some volumetric instability
is still visible during the evolution, and reducing the ratio further to
$h_s/h_x\sim\tfrac{1}{2}$ instead produced volumetric locking, with the ring
failing to stretch at all. The pressure-loaded inflating ring thus remains a
work in progress: the interplay between the solid-to-fluid mesh ratio, the
non-matching quadrature and the volumetric stabilization near a pressure
discontinuity needs to be analysed more thoroughly, and we return to this in
Section~\ref{sec:conclusions}.

\section{Conclusions}
\label{sec:conclusions}

We have presented a volumetrically stabilized mixed formulation of the
distributed-Lagrange-multiplier finite element immersed boundary method for
fluid--structure interaction problems in which a fully incompressible
hyperelastic solid is immersed in an incompressible Newtonian fluid. Starting
from the fictitious-domain FE-IBM of Boffi et al.~\cite{boffilm}, we observed
that a direct discretization of the deviatoric solid stress fails to enforce the
Lagrangian incompressibility constraint $J=1$ and produces spurious volumetric
instabilities and locking. To remedy this, we augmented the immersed solid stress
with a volumetric contribution derived from a dilatational strain energy and
introduced a solid pressure field, enforced weakly, that acts as the Lagrange
multiplier of the incompressibility constraint within the solid. In contrast to
the finite-difference IFED stabilization of~\cite{boyce1}, the solid pressure
here is a genuine mixed finite element unknown, posed and discretized entirely
within the fully variational immersed formulation. We derived the mixed
formulation for a fully incompressible neo-Hookean material, gave its
unconditionally stable semi-implicit time discretization and its finite element
space discretization, and described a reusable implementation as an
\code{ImmersedBoundary} physics kernel in the \grins multiphysics framework,
including the overlapping fluid--solid map and the parallel coupling
infrastructure required by the non-matching discretization.

The formulation was verified on three FSI benchmarks. For an elliptically
displaced thick ring the mixed formulation removes the volumetric collapse
exhibited by the unstabilized method and recovers the circular equilibrium at
every refinement level. For a radially stretched incompressible ring it sustains
the physically expected pressure difference across the ring wall, provided the
solid-to-fluid mesh ratio is sufficiently small and the mesh sufficiently
refined. For a disk falling under gravity it reproduces the empirical terminal
velocity of a rigid disk to within about~$1\%$ across a range of solid densities
and fluid viscosities, while requiring no remeshing---in contrast to a
body-fitted ALE treatment of the same problem. Throughout, we found that a
moderate amount of extra quadrature is needed to integrate the non-matching
fluid--solid coupling accurately, consistent with quadrature-error analyses of
the fictitious-domain FE-IBM.

Two limitations point to clear directions for future work. First, the accuracy of
the stabilized incompressibility enforcement depends on the solid-to-fluid
mesh-size ratio $h_s/h_x$, and a pressure-loaded inflating-ring test exposed a
regime in which the method either fails to maintain incompressibility (large
$h_s/h_x$) or locks (small $h_s/h_x$). A more thorough analysis of the interplay
between $h_s/h_x$, the non-matching quadrature, and the volumetric stabilization
near pressure discontinuities is needed; recent alternatives such as
divergence-free composite B-spline kernels~\cite{griffithcbs} and the ongoing
theoretical developments of the fictitious-domain FE-IBM~\cite{alshehri2025}
offer promising avenues. Second, we plan to validate the method against the
established Turek--Hron FSI benchmark~\cite{turek}, to extend the numerical
experiments to three dimensions, and to address co-dimension-one structures and
very large deformations, with the long-term goal of simulating the
fluid--structure interaction of a parachute deployment system.

\section*{Acknowledgements}
The authors gratefully acknowledge the members of the Applied Computational
Engineering and Science (ACES) Laboratory at the University at Buffalo for many
helpful discussions, and the developers of the \grins and \libmesh libraries.

\clearpage 

\bibliographystyle{elsarticle-num}
\bibliography{refs}

\begin{thebibliography}{10}
\expandafter\ifx\csname url\endcsname\relax
  \def\url#1{\texttt{#1}}\fi
\expandafter\ifx\csname urlprefix\endcsname\relax\def\urlprefix{URL }\fi
\expandafter\ifx\csname href\endcsname\relax
  \def\href#1#2{#2} \def\path#1{#1}\fi

\bibitem{ale}
T.~J.~R. Hughes, W.~K. Liu, T.~K. Zimmermann, {L}agrangian--{E}ulerian finite
  element formulation for incompressible viscous flows, Computer Methods in
  Applied Mechanics and Engineering 29~(3) (1981) 329--349.

\bibitem{hyperelastic}
D.~Boffi, L.~Gastaldi, L.~Heltai, C.~S. Peskin, On the hyper-elastic
  formulation of the immersed boundary method, Computer Methods in Applied
  Mechanics and Engineering 197~(25--28) (2008) 2210--2231.

\bibitem{peskin}
C.~S. Peskin, The immersed boundary method, Acta Numerica 11 (2002) 479--517.

\bibitem{heltaithesis}
L.~Heltai, The finite element immersed boundary method, Ph.D. thesis,
  Universit\`a di Pavia, Dipartimento di Matematica (2006).

\bibitem{boffi2003}
D.~Boffi, L.~Gastaldi, A finite element approach for the immersed boundary
  method, Computers \& Structures 81~(8--11) (2003) 491--501.

\bibitem{boffi2011}
D.~Boffi, N.~Cavallini, L.~Gastaldi, The finite element immersed boundary
  method with different fluid and solid densities, Mathematical Models and
  Methods in Applied Sciences 21~(12) (2011) 2523--2550.

\bibitem{boffilm}
D.~Boffi, N.~Cavallini, L.~Gastaldi, The finite element immersed boundary
  method with distributed {L}agrange multiplier, SIAM Journal on Numerical
  Analysis 53~(6) (2015) 2584--2604.

\bibitem{glow1}
V.~Girault, R.~Glowinski, T.~W. Pan, A fictitious domain method with
  distributed {L}agrange multipliers for the {S}tokes problem, Applied
  Nonlinear Analysis (Kluwer/Plenum) (1999) 159--174.

\bibitem{glow2}
R.~Glowinski, Y.~Kuznetsov, Distributed {L}agrange multipliers based on
  fictitious domain method for second order elliptic problems, Computer Methods
  in Applied Mechanics and Engineering 196~(8) (2007) 1498--1506.

\bibitem{alshehri2025}
N.~Alshehri, D.~Boffi, F.~Credali, L.~Gastaldi, Advances on finite element
  discretization of fluid--structure interaction problems, Arabian Journal of
  Mathematics 15~(1) (2025) 27--51.

\bibitem{grins}
P.~T. Bauman, R.~H. Stogner, {GRINS}: A multiphysics framework based on the
  {libMesh} finite element library, SIAM Journal on Scientific Computing 38~(5)
  (2016) S78--S100.

\bibitem{libmesh}
B.~S. Kirk, J.~W. Peterson, R.~H. Stogner, G.~F. Carey, {libMesh}: A {C}++
  library for parallel adaptive mesh refinement/coarsening simulations,
  Engineering with Computers 22~(3--4) (2006) 237--254.

\bibitem{boyce1}
B.~Vadala-Roth, S.~Rossi, S.~Prakash, N.~A. Barrett, B.~E. Griffith,
  Stabilization approaches for the hyperelastic immersed boundary method for
  problems of large-deformation incompressible elasticity, Computer Methods in
  Applied Mechanics and Engineering 365 (2020) 112978.

\bibitem{griffithcbs}
L.~Li, C.~Gruninger, J.~H. Lee, B.~E. Griffith, Local divergence-free immersed
  finite element--difference method using composite {B}-splines, arXiv preprint
  arXiv:2412.15408.

\bibitem{stein}
U.~Brink, E.~Stein, On some mixed finite element methods for incompressible and
  nearly incompressible finite elasticity, Computational Mechanics 19~(1)
  (1996) 105--119.

\bibitem{rossi}
E.~Karabelas, G.~Haase, G.~Plank, C.~M. Augustin, Versatile stabilized finite
  element formulations for nearly and fully incompressible solid mechanics,
  Computational Mechanics 65~(1).

\bibitem{gunz}
M.~D. Gunzburger, Finite Element Methods for Viscous Incompressible Flows: A
  Guide to Theory, Practice, and Algorithms, Academic Press, 1989.

\bibitem{cmbook}
O.~Gonzalez, A.~M. Stuart, A First Course in Continuum Mechanics, Cambridge
  University Press, 2008.

\bibitem{heltaiibm}
L.~Heltai, F.~Costanzo, Variational implementation of immersed finite element
  methods, Computer Methods in Applied Mechanics and Engineering 229--232
  (2012) 110--127.

\bibitem{roy}
S.~Roy, L.~Heltai, F.~Costanzo, Benchmarking the immersed finite element method
  for fluid--structure interaction problems, Computers \& Mathematics with
  Applications 69~(10) (2015) 1167--1188.

\bibitem{zhang}
L.~T. Zhang, M.~Gay, Immersed finite element method for fluid--structure
  interactions, Journal of Fluids and Structures 23~(6) (2007) 839--857.

\bibitem{clift}
R.~Clift, J.~R. Grace, M.~E. Weber, Bubbles, Drops, and Particles, Academic
  Press, New York, 1978.

\bibitem{turek}
S.~Turek, J.~Hron, Proposal for numerical benchmarking of fluid--structure
  interaction between an elastic object and laminar incompressible flow, in:
  Fluid--Structure Interaction, Vol.~53 of Lecture Notes in Computational
  Science and Engineering, Springer, 2006, pp. 371--385.

\end{thebibliography}

\end{document}